\documentclass[a4paper]{article}
\usepackage[
  top=1in,
  bottom=1in,
  left=1in,
  right=1in
]{geometry}

\usepackage{authblk}
\usepackage{setspace}
\usepackage{array}
\usepackage{makecell}
\usepackage{amsmath, amssymb, booktabs, mathtools}
\usepackage{subfig}
\usepackage[font=small,labelfont=bf]{caption}
\usepackage{threeparttable}
\usepackage{float}
\usepackage[square, sort&compress, numbers]{natbib}
\usepackage{bm}
\usepackage{stmaryrd}
\usepackage{multirow}

\usepackage[inkscapelatex=false]{svg}
\usepackage{lipsum}

\usepackage[ruled]{algorithm2e}
\SetKwComment{Comment}{/* }{ */}
\usepackage[hang,flushmargin]{footmisc}
\makeatletter
\newcommand{\algorithmfootnote}[2][\footnotesize]{%
  \let\old@algocf@finish\@algocf@finish  
  \def\@algocf@finish{\old@algocf@finish  
    \leavevmode\rlap{\begin{minipage}{\linewidth}
    #1#2
    \end{minipage}}%
  }%
}

\usepackage{graphicx} 
\usepackage{hyperref} 
\usepackage[nameinlink,capitalize]{cleveref}
\hypersetup{
  colorlinks   = true,
  urlcolor     = blue,
  linkcolor    = blue,
  citecolor    = blue
}
\newcommand{\correspondingA}{*}

\title{Homotopy continuation of viscoelastic waveguide dispersion curves: from intra‑manifold tracking to inter‑manifold transport}

\author[]{Dong Xiao\textsuperscript{*}\href{https://orcid.org/0009-0006-7609-7832}{\includegraphics[scale=0.04]{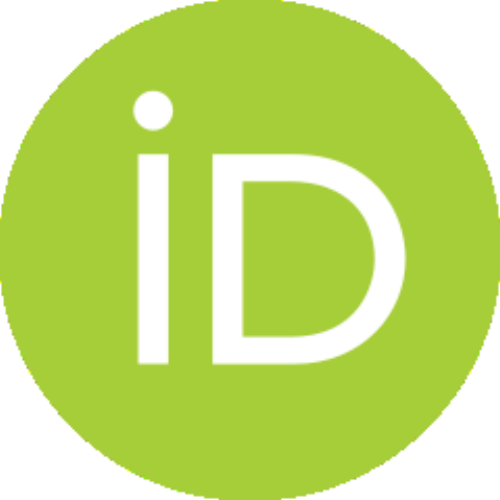}}}
\author[]{Zahra Sharif Khodaei\href{https://orcid.org/0000-0001-5106-2197}{\includegraphics[scale=0.04]{orcidicon.pdf}}}
\author[]{M. H. Aliabadi\href{https://orcid.org/0000-0002-2883-2461}{\includegraphics[scale=0.04]{orcidicon.pdf}} }

\affil[]{\normalsize \slshape  Department of Aeronautics, Imperial College London, South Kensington, London SW7 2AZ, United Kingdom.}
\date{\vspace*{-2\baselineskip}}  
\date{\vspace*{-1.5cm}}

\begin{document}
\maketitle
\footnotetext[0]{ \textsuperscript{\correspondingA}Corresponding author}
\footnotetext[1]{Email addresses: d.xiao21@imperial.ac.uk (D. Xiao); z.sharif-khodaei@imperial.ac.uk (Z. Sharif-Khodaei); m.h.aliabadi@imperial.ac.uk (M.H. Aliabadi)}

\renewcommand{\abstractname}{Abstract}
\begin{abstract}{\normalsize \onehalfspacing
Conventional mode tracking operates in the dark: it traces dispersion branches on the non‑Hermitian eigenvalue manifold using only local continuity, unaware of the global Riemann‑sheet topology. When exceptional points (EPs) lie close to the real frequency axis, the eigenvector similarity that local trackers rely on degrades, and mode tracking becomes unreliable, failing silently. This paper replaces blind intra‑manifold tracking with inter‑manifold transport. A material attenuation parameter s in [0,1] continuously maps the target lossy problem to an auxiliary lossless one whose Hermitian eigenvalue problem yields a well‑posed anchor manifold on which each dispersion branch possesses a globally unique and continuous identity. These identities are defined once on the elastic anchor and then transported to the viscoelastic target via predictor–corrector homotopy continuation; as long as the path avoids all EPs, branch identity is preserved throughout the transport. For any mode pair whose EPs have not crossed the real frequency axis (Type I), the transported identities are inherited automatically. In contrast, when an EP crosses the real axis and becomes Type II, the topology differs from the elastic anchor and a label swap is required. The framework is validated on symmetric and unsymmetric laminates, with most cases at loss factors of 0.003 to 0.02; for all Type~I pairs in these cases the identities are inherited without alteration. For a challenging unsymmetric laminate at 0.05, several EP pairs have become Type~II, yet the homotopy transport still produces numerically accurate solutions. Two diagnostic signatures—an extremely sharp imaginary‑part crossing and a marked discrepancy between spectral group velocity and energy flux velocity—identify where the underlying EP topology demands a label swap.
}
\end{abstract}

\renewcommand{\abstractname}{Keywords}
\begin{abstract}{Viscoelastic waveguides; Non-Hermitian eigenvalue problems; Homotopy continuation; Exceptional points; Type I/II EP topology; Inter‑manifold transport.}
\end{abstract}

\newpage

\section{Introduction}
Ultrasonic guided waves have emerged as a powerful tool for the non‑destructive evaluation (NDE) and structural health monitoring (SHM) of elongated engineering structures, owing to their ability to propagate over long distances and their sensitivity to a variety of defects~\cite{rose_ultrasonic_2014, giurgiutiu_structural_2014}. The successful application of guided‑wave‑based techniques relies on a thorough understanding of dispersion characteristics—the relationship between frequency, phase velocity, group velocity, and, for lossy materials, attenuation~\cite{auld_acoustic_1990}. For waveguides with simple geometries and elastic material behavior, analytical or semi‑analytical solutions are available~\cite{meleshko_elastic_2009}. 
Modern engineering structures, however, increasingly employ materials with inherent viscoelastic properties, such as carbon fiber‑reinforced polymer (CFRP) with lossy matrices, as well as adhesive layers and damping treatments that exhibit viscoelastic behavior even when the base material is elastic~\cite{nashif_vibration_1985, hosten_transfer_1993}. Consequently, accurate dispersion analysis of viscoelastic waveguides of arbitrary cross‑section has become essential for NDE/SHM of laminated composites, polymeric pipes, and structures with damping.

A wide variety of methods have been developed for this purpose. Matrix methods based on partial‑wave decomposition—including the Global Matrix Method (GMM), Transfer Matrix Method (TMM), and Stiffness Matrix Method (SMM) and its variants—enforce continuity conditions at layer interfaces and have been implemented in widely used software packages~\cite{lowe_matrix_1995, maghsoodi_calculation_2014, wang_stable_2001, rokhlin_stable_2002, pavlakovic_disperse_1997, huber_dispersion_2024, huber_stiffness_2024}. The Semi‑Analytical Finite Element (SAFE) method discretizes only the cross‑section, reducing the problem to an eigenvalue problem for complex wavenumbers at each frequency~\cite{gavric_finite_1994, castaings_finite_2008, rose_semi-analytical_2014}, and has been extended to damped, piezoelectric, prestressed, thin‑walled orthotropic, and multi‑material waveguides, as well as integrated into user‑friendly software tools~\cite{bartoli_modeling_2006, marzani_semi-analytical_2008, mu_guided_2008, mazzotti_coupled_2013, cortes_dispersion_2008, mazzotti_guided_2012, mazzotti_wave_2016, bocchini_graphical_2011, liu_modelling_2024}. Alternative global discretization techniques such as the spectral collocation method (SCM) and the Legendre Polynomial Method (LPM) offer spectral accuracy using global basis functions and have also been applied to viscoelastic media~\cite{adamou_spectral_2004, quintanilla_guided_2015, dahmen_investigation_2016}.  
Recent comparative studies confirm that these methods exhibit distinct trade‑offs in accuracy, efficiency, and robustness for viscoelastic waveguides. Orta et al.~\cite{orta_comparative_2022} benchmarked six methods on multi‑layer laminates, highlighting that while matrix methods and SAFE deliver comparable accuracy, their computational costs differ markedly; Quiroga et al.~\cite{quiroga_evaluation_2025} further demonstrated that SAFE and LPM achieve sub‑1\% errors for fundamental modes up to 100 kHz while significantly outperforming GMM in speed. Barazanchy and Giurgiutiu~\cite{barazanchy_comparative_2017} provided a concise user‑oriented overview of unified analytic, SAFE, and DISPERSE methods for composites.  

Despite this methodological diversity, all existing approaches share two intertwined difficulties when applied to viscoelastic waveguides: (i) solving a non‑Hermitian eigenvalue problem at each frequency, and (ii) reliably tracing the resulting eigenmodes across frequencies to form continuous dispersion curves (see also \cref{app:challenges} for a comprehensive discussion).

The first difficulty is the solution of the non‑Hermitian eigenvalue problem itself. Material damping introduces complex arithmetic and destroys the Hermitian structure of the system matrices, leading to a polynomial eigenvalue problem in the complex wavenumber. The various discretization strategies—SAFE with linearization, SCM with a companion matrix, or matrix methods—are all ``formulations'' that construct this eigenvalue problem; they must be paired with a numerical solver. Several solver classes exist, each with intrinsic limitations. Dense direct eigensolvers (e.g., \texttt{eig} in MATLAB) compute all eigenvalues at once, but their $\mathcal{O}(n^3)$ cost and $\mathcal{O}(n^2)$ memory become prohibitive for large finite‑element models. Iterative eigensolvers such as shift‑invert Arnoldi mitigate this scaling, yet they introduce heuristic parameter selection: the choice of shift values critically determines which modes are found, and multiple shifts are often required to explore the complex plane, risking omission of relevant solutions~\cite{lehoucq_deflation_1996, tisseur_quadratic_2001}. Contour‑integral methods~\cite{sakurai_projection_2003} can robustly extract all eigenvalues within a specified region, but at a higher computational cost. Complex root‑searching schemes~\cite{orta_comparative_2022, quiroga_evaluation_2025}, while accurate for isolated modes, become increasingly cumbersome in dense spectral regions or when modes veer and cross.

The second difficulty, mode tracking, is equally critical. In viscoelastic waveguides, the non‑Hermitian character complicates mode veering and degeneracy: eigenvectors exchange character rapidly over narrow frequency intervals, challenging even sophisticated tracking algorithms, let alone tools relying on the Modal Assurance Criterion (MAC)~\cite{seyranian_multiparameter_2003, allemang_modal_2003}. Frequency‑continuation methods have been used to solve nonlinear eigenvalue problems and to trace modal branches in structural dynamics~\cite{akoussan_numerical_2018, liu_high-efficient_2024, boudaoud_numerical_2009, ziapkoff_high_2024}, and have also been applied to dispersion mode tracking~\cite{allgower_numerical_1990, maruyama_continuation_2025}. However, their direct application at the target viscoelastic state introduces fundamental difficulties: they require a high‑quality starting solution on the non‑Hermitian manifold and are highly sensitive to exceptional points (EPs), where eigenvalues coalesce and the continuation Jacobian becomes singular. Thus, direct intra‑manifold continuation on the target viscoelastic state is fragile and often fails silently. 

Herein lies the fundamental structural limitation: all conventional methods, regardless of their formulation or solver, fix the material state at the target viscoelastic configuration and attempt to track modes ``on'' the same non‑Hermitian eigenvalue manifold—that is, they perform ``intra‑manifold tracking''. They rely on local information such as eigenvector similarity, tangent prediction, or correlation metrics, and possess no awareness of the global Riemann‑sheet topology on which the eigenvalues reside. As exceptional points migrate toward the real frequency axis, the eigenvectors of interacting modes collapse into nearly collinear directions. Correlation metrics such as the MAC—which rely on vector inner products—accordingly lose discriminative power, and local trackers may silently exchange branch identities, producing numerically smooth curves that conceal physically incorrect labels. In effect, conventional mode tracking operates in the dark. 

This topological ignorance is the root cause of fragility, but the observable behavior of dispersion curves is not accidental. Crossing theory for non‑Hermitian systems~\cite{keck_unfolding_2003} shows that a pair of interacting EPs can adopt two distinct configurations relative to the real frequency axis. In a Type~I configuration, the EPs lie on opposite sides of the axis; along the real axis the real parts of the eigenfrequencies veer while the imaginary parts cross. In a Type~II configuration, both EPs lie on the same side of the axis, leading to a crossing of real parts and veering of imaginary parts. These veering and crossing patterns are direct projections of the underlying EP topology~\cite{heiss_phases_1999,heiss_exceptional_2004,ghienne_beyond_2020}. For the lossless elastic waveguide, every mode-veering event is the spectral footprint of a Type~I EP pair symmetric about the real axis. As material damping increases, EPs migrate; if an EP crosses the real axis, a Type~II configuration emerges. This continuous dependence on the damping level—and the one‑to‑one correspondence between EP configuration and observable crossing type—constitutes the physical foundation for the homotopy framework we propose.

To circumvent the topological blindness that undermines existing methods, this work adopts an ``inter-manifold transport'' strategy. Rather than tracing modes directly on the lossy non-Hermitian target manifold, we first establish each modal identity on the Hermitian manifold of the lossless elastic waveguide—a well-posed ``anchor manifold'' where every dispersion branch possesses a globally unique and continuous identity \cite{xiao_rigorous_2026, gravenkamp_notes_2023}. These pre-established labels are then transported to the viscoelastic target along a material homotopy path. 

The transport is realized by a linear material homotopy, $\mathbf{C}(s) = \mathbf{C}' + is\mathbf{C}''$, in which the scalar attenuation factor \(s \in [0,1]\) interpolates between the elastic anchor (\(s=0\)) and the target viscoelastic state (\(s=1\)). Identities are defined at \(s=0\)  and propagated via predictor–corrector homotopy continuation. As long as the path avoids all exceptional points, branch identity is preserved—a guarantee rooted in analytic perturbation theory \cite{kato_perturbation_1995}. Whether the transported identities are inherited correctly at the target state (\(s=1\)) depends on the EP topology. For mode pairs governed by Type I EPs (EPs on opposite sides of the real frequency axis), the target manifold is compatible with the anchor, and the identities are inherited automatically. In contrast, when an EP crosses the real frequency axis and becomes Type II, the topology has changed relative to the anchor, and a label swap is required to restore physically correct mode identities.

This architecture decouples what was previously a single tightly coupled challenge—solving a non‑Hermitian eigenvalue problem while simultaneously tracking modes—into two sequential, well‑posed tasks: reliable mode identification in the Hermitian regime, followed by safe identity transport along a topologically informed path. The first task—robustly assigning globally unique identities to every dispersion branch on the lossless elastic manifold—has been rigorously addressed in our prior work \cite{xiao_rigorous_2026}, where an adaptive wavenumber sampling algorithm grounded in the Wigner–von Neumann non‑crossing rule provides provable guarantees for correct mode tracking in Hermitian SAFE formulations. The second task—transporting these identities to the viscoelastic target—is accomplished by the inter‑manifold homotopy continuation and topological inheritance rules described above; the analytic structure of the homotopy path and the EP classification (Type I versus Type II) govern whether labels are retained automatically or must be swapped post hoc, as detailed in \cref{sec:theory_guarantee}. 

Because the homotopy is performed at fixed frequency, any frequency-dependent damping model reduces to a scaled frequency-independent (hysteretic) damping at each step; the framework therefore accommodates both hysteretic and frequency-dependent damping models under a single unified parameterization, and applies without modification to waveguides of arbitrary cross-section discretizable by the SAFE method. Owing to the lack of reference solutions for frequency-dependent viscoelastic models, the present paper adopts a hysteretic damping model and validates the framework numerically on symmetric laminates, unsymmetric laminates, and an L-shaped bar.

The contributions of this work are threefold:
\begin{enumerate}
\item \textbf{Paradigm shift: from intra-manifold tracking to inter-manifold transport.}
We abandon the conventional approach of tracking modes on the non-Hermitian target manifold using only local similarity and instead define mode identities on the well-posed Hermitian anchor manifold, then transport them to the viscoelastic target along a safe homotopy path. This decouples a simultaneous, tightly coupled challenge into two sequential, well-posed tasks.

\item \textbf{Physical criterion for mode-label inheritance via non-Hermitian EP topology.}
By synthesizing non-Hermitian EP physics and crossing theory with guided-wave dispersion analysis, we introduce a physically grounded criterion for label inheritance after homotopy transport. The Type I/II EP topology—governing whether interacting EPs lie on opposite or the same side of the real frequency axis—directly dictates the observable crossing or veering patterns of dispersion branches. This correspondence yields a simple inheritance rule: Type I configurations retain anchor labels automatically, whereas Type II configurations require a label swap. Because explicit EP tracking in the complex frequency plane of a non-Hermitian waveguide demands dedicated algorithms beyond the scope of this work, the Type II transition cannot be predicted a priori here; we therefore provide two empirical diagnostic signatures—an extremely sharp imaginary-part crossing and a marked spectral–energetic velocity discrepancy—to flag when the topology has transitioned to Type II.

\item \textbf{Scalable and parallelizable computational framework.}
We develop a predictor–corrector homotopy continuation algorithm with adaptive step-size control that propagates sparsely selected key solutions from the Hermitian anchor to the viscoelastic target. The framework is inherently parallelizable, applies without modification to arbitrary cross-sections discretizable by SAFE, and extends—under the analyticity condition established herein—to frequency-dependent damping models by treating them as scaled hysteretic damping at each fixed-frequency step.
\end{enumerate}

The remainder of this paper is organized as follows. \Cref{sec:methodology} presents the mathematical formulation of the SAFE method for viscoelastic waveguides with hysteretic damping, details the construction of the material homotopy, describes its numerical implementation. It also discusses the theoretical foundations of branch identity continuity along the homotopy path, as well as the correspondence of EP topology to crossing type in non-hermitian system. \Cref{sec:validation} presents numerical examples that validate the accuracy and efficiency of the proposed method. \Cref{sec:discussion} discusses the computational efficiency, scalability, robustness, and advantages over existing methods, including a dedicated analysis of the effect of loss factor magnitude and eigengap on performance. Finally, \Cref{sec:conclusion} concludes the paper and outlines limitations and future work.

\section{Methodology}
\label{sec:methodology}
\subsection{SAFE formulation for viscoelastic waveguides: from Hermitian to non-Hermitian}
The Semi-Analytical Finite Element (SAFE) method provides an efficient computational framework for calculating dispersion characteristics in waveguides of uniform cross-section \cite{bartoli_modeling_2006, marzani_semi-analytical_2008, rose_semi-analytical_2014}. Following the formulation detailed in our previous work \cite{xiao_rigorous_2026}, the displacement field within an element is expressed as a harmonic function along the propagation direction:
\begin{equation}
    \begin{aligned}
        \mathbf{u}^{(e)}(x,y,z,t) = \mathbf{N}^{(e)}(y, z) \mathbf{q}^{(e)} e^{i(kx - \omega t)}, 
    \end{aligned}
\end{equation}
where $\mathbf{N}^{(e)}(y,z)$ contains the shape functions, $\mathbf{q}^{(e)}$ is the vector of nodal displacements, $k$ and $\omega$ denote the wavenumber and angular frequency, respectively, $i$ is the imaginary unit, and $t$ is time. Substituting this representation into the elastodynamic equations and applying the principle of virtual work yields the quadratic eigenvalue problem in the wavenumber $k$:
\begin{equation} \label{eq:SAFE_quad_eigen}
    \begin{aligned}
       \mathbf{D}(k, \omega)\mathbf{q} = (\mathbf{K}_1 + ik\mathbf{K}_2+ k^2\mathbf{K}_3 -\omega^2 \mathbf{M})\mathbf{q} =  \mathbf{0},
    \end{aligned}
\end{equation}
where $\mathbf{D}(k, \omega) =  (\mathbf{K}_1 + ik\mathbf{K}_2+ k^2\mathbf{K}_3 -\omega^2 \mathbf{M})$ is dynamic stiffness matrix of SAFE system, $\mathbf{q}$ collects the nodal degrees of freedom of the entire cross-sectional discretization. The global matrices are assembled from element contributions:
\begin{equation} \label{eq:K_j_int}
    \begin{aligned}
    \mathbf{K}_1 &= \int_A \mathbf{B}_0^T\mathbf{C}\mathbf{B}_0 dA, \\
    \mathbf{K}_2  &= \int_A \left(\mathbf{B}_0^T\mathbf{C}\mathbf{B}_1 - \mathbf{B}_1^T\mathbf{C}\mathbf{B}_0\right) dA, \\
    \mathbf{K}_3  &= \int_A \mathbf{B}_1^T\mathbf{C}\mathbf{B}_1 dA, \\
    \mathbf{M}  &= \int_A \rho\mathbf{N}^T\mathbf{N} dA.
    \end{aligned}
\end{equation}
Here $\mathbf{B}_0$ and $\mathbf{B}_1$ are strain–displacement matrices associated with the zeroth and first derivatives with respect to $x$, $\mathbf{C}$ is the material stiffness tensor, $\mathbf{N}$ is the shape function matrix, $\rho$ is the mass density, and $A$ denotes the cross-sectional area. For a lossless (elastic) waveguide, $\mathbf{C}$ is real and symmetric, rendering $\mathbf{K}_1$, $\mathbf{K}_3$ and $\mathbf{M}$ symmetric positive-definite, while $\mathbf{K}_2$ is skew‑symmetric. Consequently, for any real wavenumber $k$, the stiffness matrix $\mathbf{K}(k) = \mathbf{K}_1 + i k \mathbf{K}_2 + k^2 \mathbf{K}_3$ is complex Hermitian ($\mathbf{K}(k)=\mathbf{K}^{\dag}(k)$), leading to a well-behaved eigenvalue problem.

The present work extends this formulation to account for viscoelastic behavior. For the hysteretic (frequency-independent) damping model \cite{maia_reflections_2009}, the material stiffness tensor is complex:
\begin{equation} \label{eq:hysteretic}
    \begin{aligned}
    \mathbf{C} = \mathbf{C}' + i \mathbf{C}''.
    \end{aligned}
\end{equation}
Substituting \cref{eq:hysteretic} into \cref{eq:K_j_int} yields complex stiffness matrices:
\begin{equation} \label{eq:K_j_comp}
    \begin{aligned}
    \mathbf{K}_j = \mathbf{K}'_j + i \mathbf{K}''_j, \; j \in \{1,2,3\},
    \end{aligned}
\end{equation}
where $\mathbf{K}_j'$ are assembled from $\mathbf{C}'$, and $\mathbf{K}_j''$ are assembled once from $\mathbf{C}''$. In this case, the dynamic stiffness matrix $\mathbf{K}(k)$ is no longer Hermitian. The resulting non-Hermitian eigenvalue problem requires tracking eigenvalues and eigenvectors across a Riemann surface whose topology is governed by exceptional points (EPs). Key concepts from non-Hermitian dispersion analysis—including the Riemann surface, branch identity, EPs, dispersion curves, mode tracking and the distinction between spectral group velocity and energy flux velocity—are summarised in \cref{app:terminology}. The adaptive homotopy framework presented next decouples this challenging non-Hermitian problem by shifting mode tracking to the elastic limit, where the system is Hermitian and modal identities are unambiguous.

\subsection{Adaptive homotopy continuation for non-Hermitian dispersion analysis: inter-manifold transport with theoretical foundations}
The proposed framework adopts an inter-manifold transport strategy: instead of tracing modes directly on the lossy non-Hermitian target manifold, we establish globally unique modal identities on the Hermitian elastic anchor ($s=0$) and transport them to the viscoelastic target ($s=1$) along a material homotopy path. This decouples the coupled challenge of non-Hermitian eigensolving and mode tracking into two sequential tasks—robust mode identification in the lossless regime, followed by safe identity propagation via predictor–corrector continuation in the attenuation parameter $s \in [0,1]$.

Whether the transported identities are inherited correctly at $s=1$ is governed by the EP topology. Under a Type I configuration (EPs on opposite sides of the real frequency axis), the target manifold is compatible with the anchor, and physical labels are retained automatically. Under a Type II configuration (an EP has crossed the real axis), the topology has changed relative to the anchor, and a post-hoc label swap is required. Provided the homotopy path avoids all EPs in the $s$-plane, analytic perturbation theory guarantees branch identity continuity and a one-to-one correspondence between elastic and viscoelastic solutions.

\begin{figure}[htb]
\centering
\includegraphics[width=1.02\columnwidth]{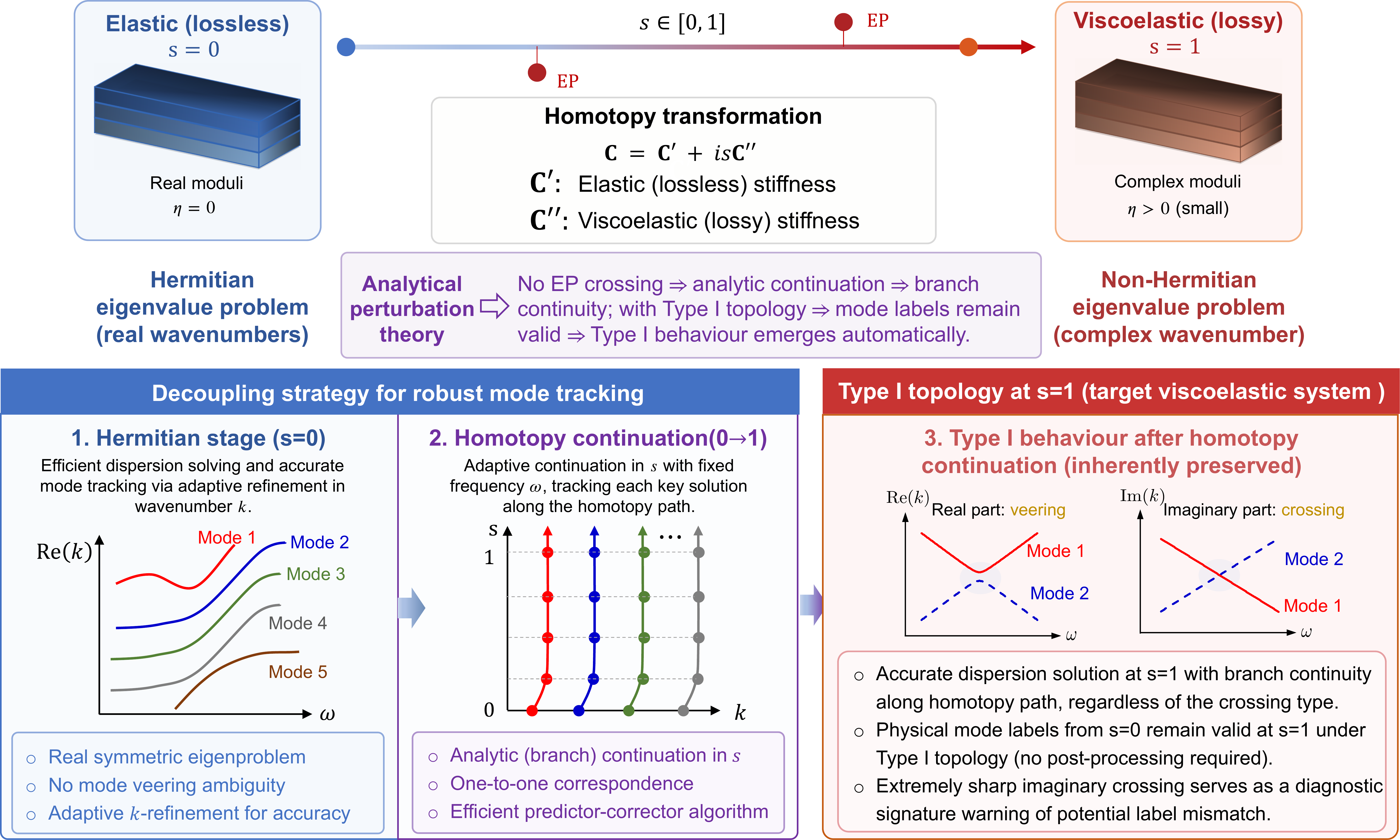}
\caption{Overview of the proposed adaptive homotopy continuation framework. A material homotopy parameterised by $s \in [0,1]$ continuously maps the elastic (lossless, Hermitian) system at $s=0$ to the viscoelastic (lossy, non-Hermitian) system at $s=1$. Mode identification is performed in the Hermitian setting and propagated to the viscoelastic regime via homotopy continuation. When the target system exhibits a Type I EP topology, the physical mode labels remain valid at $s=1$, and the resulting dispersion curves naturally exhibit Type I behaviour without post-processing.}
\label{fig:method_overview}
\end{figure}

\Cref{fig:method_overview} schematizes the overall strategy. Stage 1 performs Hermitian mode tracking at $s=0$ , yielding a complete set of dispersion branches with unambiguous physical labels. These serve as starting points for Stage 2, where predictor–corrector continuation traces each branch along the homotopy path $s=0 \to s=1$  at fixed frequency $\omega$ , allowing the complex wavenumber $k$ and eigenvector $\mathbf{q}$ to adapt to the evolving material properties. For systems retaining a Type I EP topology, the labels established at the elastic anchor remain valid at $s=1$, and the characteristic Type I behaviour—real-part veering with imaginary-part crossing-emerges automatically without post-processing. The detailed homotopy construction, lossless-stage algorithm, and predictor-corrector implementation are given in \cref{sec:homotopy_construction,sec:lossless_stage,sec:material_homotopy}.

\subsubsection{Construction of the material homotopy}
\label{sec:homotopy_construction}
For the hysteretic model, the homotopy is constructed by introducing a parameter $s \in [0,1]$ that scales the material loss:
\begin{equation} \label{eq:homotopy_C}
    \begin{aligned}
    \mathbf{C}(s) = \mathbf{C}' + i s  {\mathbf{C}}''.
    \end{aligned}
\end{equation}
This linear interpolation preserves the physical interpretation of $s$ as a "loss factor" and ensures that at $s=0$, the material is purely elastic with real stiffness $\mathbf{C}'$, while at $s=1$, the full viscoelastic behavior $\mathbf{C}' + i {\mathbf{C}}''$ is recovered. Substituting \cref{eq:homotopy_C} into the \cref{eq:K_j_int} yields parameter-dependent stiffness matrices $\mathbf{K}_1(s)$, $\mathbf{K}_2(s)$, and $\mathbf{K}_3(s)$, all linear in $s$ due to the linearity of the integral operators:
\begin{equation} \label{eq:Kj_t_linear}
    \begin{aligned}
    \mathbf{K}_j(s) = \mathbf{K}_j' + i s {\mathbf{K}}_j'', \; j \in {1,2,3}.
    \end{aligned}
\end{equation}
The partial derivative with respect to the homotopy parameter follows directly:
\begin{equation} \label{eq:Kj_deriv}
    \begin{aligned}  
    \frac{\partial \mathbf{K}_j}{\partial s} = i {\mathbf{K}}_j'', \; j \in {1,2,3}.
    \end{aligned}
\end{equation}
which is a constant matrix, independent of the frequency $\omega$ and homotopy parameter $s$. This linear dependence on $s$ ensures that the Jacobian matrix $\partial \mathbf{G}/\partial s$ (see \cref{eq:G_t}) can be evaluated efficiently without recomputing element integrals during path tracking.

For a fixed frequency $\omega$, the homotopy defines a family of nonlinear eigenvalue problems:
\begin{equation}\label{eq:homotopy_F}
    \begin{aligned}
    \mathbf{F}(k, \mathbf{q}, s) = [\mathbf{K}_1(s) + i k \mathbf{K}_2(s) + k^2 \mathbf{K}_3(s) - \omega^2 \mathbf{M}] \mathbf{q} = \mathbf{0}, s \in [0, 1].
    \end{aligned}
\end{equation}
At $s=0$, the system is Hermitian; at $s=1$, it is the target non-Hermitian system. To resolve the scaling and phase ambiguity of the eigenvector, we append a complex normalization condition:
\begin{equation} \label{eq:normalization_condition}
    \begin{aligned}
   \mathbf{q}^{\dagger}_{\text{ref}} \mathbf{q} -1=0, 
    \end{aligned}
\end{equation}
where $\mathbf{q}_{\text{ref}}$ is a fixed reference vector in the current prediction‑correction step (typically taken as the previously converged solution). This constraint fixes the projection of $\mathbf{q}$ onto the direction of $\mathbf{q}_{\text{ref}}$, thereby determining both the norm and the phase of the eigenvector uniquely. Unlike a norm‑plus‑phase constraint, which involves the non‑analytic term $\mathbf{q}^{\dagger}\mathbf{q} -1=0$, the present formulation is analytic in $\mathbf{q}$. Together with $\mathbf{F}(k, \mathbf{q}, s)$, which comprises $n$ complex equations, the extended system $\mathbf{G}(\mathbf{y}, s) = \mathbf{0}$ with $\mathbf{y} = [\mathbf{q}^T, k]^T$ forms a square system of $n+1$ complex equations in $n+1$ complex unknowns, and can be solved directly using complex Newton iterations without separating real and imaginary parts. This analyticity also simplifies the derivation of the Jacobian matrix and preserves quadratic convergence in the predictor‑corrector steps, as derived in next \cref{sec:material_homotopy}.

\subsubsection{Stage 1: hermitian solution and mode tracking}
\label{sec:lossless_stage}
The lossless stage exploits the Hermitian nature of the $s=0$ system to obtain a complete and accurately traced set of dispersion curves. For single‑parameter Hermitian eigenvalue problems, the non‑crossing rule guarantees that mode crossings and degeneracies are symmetry‑protected \cite{von_neumann_no_1929, hatton_noncrossing_1976, mead_noncrossing_1979}. Moreover, in mode veering regions there exists a finite sampling resolution sufficient to ensure reliable MAC‑based mode tracking \cite{xiao_rigorous_2026}. Building on these theoretical foundations, the authors have previously developed an adaptive mode tracking method for Hermitian SAFE systems at $s=0$ \cite{xiao_rigorous_2026}, which is adopted in the present study.

The method proceeds as follows. First, a low-resolution wavenumber sweep is performed, solving the Hermitian eigenvalue problem at discrete wavenumber points $k_p$. An error indicator $\varepsilon$ is defined based on the MAC separation between adjacent solutions \cite{xiao_rigorous_2026}:
\begin{equation} \label{eq:error_indicator}
\begin{aligned}
    \varepsilon(k_p, k_{p+1}) = 1 - \min_j \left(\max_m \left(\mathrm{MAC}[\mathbf{q}_j(k_p), \mathbf{q}_j(k_{p+1})]- \mathrm{MAC}[\mathbf{q}_j(k_p), \mathbf{q}_m(k_{p+1})]\right)\right),
\end{aligned}
\end{equation}
where $\mathrm{MAC}$ is defined in \cref{eq:MAC}. This indicator quantifies the reliability of mode tracking between successive points. In regions where $\varepsilon$ exceeds a prescribed threshold $\bar{\varepsilon}$, the wavenumber grid is refined by inserting additional sampling points. This adaptive resampling process continues iteratively until all modal branches exhibit $\varepsilon < \bar{\varepsilon}$ between consecutive points, ensuring accurate tracking even through regions of mode veering. When symmetry-protected degeneracies are detected, a subspace tracking technique based on subspace MAC is activated to maintain correct mode identification.

Wavenumber sweeping is adopted rather than frequency sweeping at the lossless stage. This choice preserves the Hermitian nature of the system: when $k$ is prescribed as a real parameter, $\mathbf{K}(k) = \mathbf{K}_1 + i k \mathbf{K}_2 + k^2 \mathbf{K}_3$ remains complex Hermitian, and the resulting eigenvalue problem yields frequencies (up to numerical precision). The eigenvectors $\mathbf{q}(k)$ are complex-valued due to the presence of the $i k \mathbf{K}_2$ term, but the Hermitian structure ensures that modes are well-separated and trackable.

The outcome of this stage is a comprehensive dataset comprising, for each mode $j$ and each sampled wavenumber $k_p$, the frequency $\omega_{j,0}(k_p)$ (real-valued up to machine precision) and the corresponding eigenvector $\mathbf{q}_{j,0}(k_p)$, with unambiguous modal connectivity established across the entire wavenumber range of interest. Critically, for each mode $j$, the condition
\begin{equation} \label{eq:sparse_condition}
\begin{aligned}
\mathrm{MAC}[\mathbf{q}_{j,0}(k_p), \mathbf{q}_{j,0}(k_{p+1})] > 1 - \bar{\varepsilon},
\end{aligned}
\end{equation}
holds for all adjacent wavenumber pairs, guaranteeing reliable mode tracking throughout the dataset.

To reduce the number of points that must be propagated through homotopy, an adaptive sparse mapping strategy is applied to the dense dataset. The goal is to select a subset of “key” solution points such that:
\begin{enumerate}
\item \textbf{Modal continuity} is preserved: the MAC between any two consecutive key points exceeds a prescribed threshold $1-\bar{\zeta}$.
\item \textbf{Geometric accuracy} of the dispersion curve is maintained: the frequency $\omega(k)$ can be accurately reconstructed via linear interpolation between key points, with an interpolation error below a specified tolerance $\bar{\gamma}$.
\end{enumerate}
A greedy thinning algorithm is employed to construct the key point subsequence. The first condition ensures reliable mode tracking; the second guarantees that the retained points faithfully represent the shape of the dispersion curve, especially in regions of high curvature or rapid frequency variation. The resulting set of key points $(\omega_{0}^{\text{key}}, k_{0}^{\text{key}}, \mathbf{q}_{0}^{\text{key}})$ is significantly sparser than the original dense dataset yet retains all information necessary for accurate representation of the dispersion curves. These key points serve as the starting points for the material homotopy path tracking.

\subsubsection{Stage 2: adaptive homotopy path tracking}
\label{sec:material_homotopy}
For each selected starting point $(\omega_{j,0}^{\text{key}}, k_{j,0}^{\text{key}}, \mathbf{q}_{j,0}^{\text{key}})$, the homotopy path $\mathbf{y}(s) = [\mathbf{q}(s)^T, k(s)]^T$ is traced from $s=0$ to $s=1$ using a predictor-corrector algorithm with arc-length parameterization. The extended complex system of size of $n+1$ is defined by combined \cref{eq:homotopy_F} and \cref{eq:normalization_condition}:
\begin{equation}
    \begin{aligned} \label{eq:G}
    \mathbf{G}(\mathbf{y}, s) = \left[\begin{matrix}
                                 \mathbf{F}(k, \mathbf{q}, s) \\
                                \mathbf{q}_{\text{ref}}^{\dagger} \mathbf{q} - 1
                                \end{matrix}\right] = \mathbf{0} \in \mathcal{C}^{(n+1)}, 
    \end{aligned}
\end{equation}
with $\mathbf{F}(k, \mathbf{q}, s)$ given by \cref{eq:homotopy_F}. The Jacobian matrix with respect to $\mathbf{y}$ is obtained by differentiating $\mathbf{G}$ with respect to $\mathbf{q}$ and $k$. Since the constraint is linear and analytic, we have:
\begin{equation}
    \begin{aligned}\label{eq:G_y}
    \frac{\partial \mathbf{G}}{\partial \mathbf{y}} = \left[\begin{matrix}
                                 \mathbf{D}(k, s) & [i\mathbf{K}_2(s) + 2k\mathbf{K}_3(s)]\mathbf{q}\\
                                \mathbf{q}_{\text{ref}}^{\dagger} & 0 \\
                                \end{matrix}\right]  \in \mathcal{C}^{(n+1) \times (n+1)},
    \end{aligned}
\end{equation}
where $\mathbf{D}(k,s) = \mathbf{K}_1(s) + i k \mathbf{K}_2(s) + k^2 \mathbf{K}_3(s) - \omega^2 \mathbf{M}$ is the dynamic stiffness matrix, which is parameterised by $k$ and $s$ ($\omega$ is a constant for each solution in homotopy continuation). The partial derivative with respect to the homotopy parameter $s$ is:
\begin{equation}
    \begin{aligned} \label{eq:G_t}
     \frac{\partial \mathbf{G}}{\partial s} = \left[\begin{matrix}
                                [\frac{\partial \mathbf{K}_1(s)}{\partial s} + ik\frac{\partial \mathbf{K}_2(s)}{\partial s}+
                                k^2\frac{\partial \mathbf{K}_3(s)}{\partial s}
                                ] \mathbf{q} \\
                                0 
                                \end{matrix}\right],
    \end{aligned}
\end{equation}
where $\partial \mathbf{K}_j/\partial s$ is given by \cref{eq:Kj_deriv}. Since ${\mathbf{K}}_j''$ are constant matrices computed during preprocessing, $\partial \mathbf{K}_j/\partial s$ can be evaluated efficiently for each frequency $\omega$ without recomputing element integrals. 

The predictor-corrector algorithm proceeds as follows. At a known point ($\mathbf{y}_p$, ${s}_{p}$), the tangent vector $\dot{\mathbf{y}}_p = \frac{\partial \mathbf{y}_p}{\partial {s}}$ is obtained by solving the linear system:
\begin{equation} \label{eq:tangent}
    \begin{aligned}
    \frac{\partial \mathbf{G}}{\partial \mathbf{y}} \dot{\mathbf{y}} +  \frac{\partial \mathbf{G}}{\partial s}  = \mathbf{0},\; \left \| \dot{\mathbf{y}} \right\|^{2}  = 1, 
    \end{aligned}
\end{equation}
which follows from differentiating $\mathbf{G}(\mathbf{y}(s), s) = 0$ with respect to $s$ and noting that the homotopy parameter $s$ is strictly increasing along the path (no turning points). A prediction is made with step size $\Delta s$:
\begin{equation} \label{eq:predictor}
    \begin{aligned}
    \mathbf{y}_{p+1}^{(0)} = \mathbf{y}_{p} + \Delta s \cdot \dot{\mathbf{y}}_{p}, \;
    {s}_{p+1} = {s}_{p} + \Delta s,
    \end{aligned}
\end{equation}
Newton iterations then correct the prediction by solving:
\begin{equation} \label{eq:newton}
    \begin{aligned}
     \frac{\partial \mathbf{G}}{\partial \mathbf{y}}(\mathbf{y}_{p+1}^{(l)},  s_{p+1}) \Delta  \mathbf{y} = - \mathbf{G}(\mathbf{y}_{p+1}^{(l)},  s_{p+1}).
    \end{aligned}
\end{equation}
updating $\mathbf{y}_{p+1}^{(l+1)} = \mathbf{y}_{p+1}^{(l)} + \Delta \mathbf{y}$ until convergence. The linear systems are solved in complex arithmetic, exploiting the analyticity of $\mathbf{G}$ to maintain quadratic convergence. This approach avoids the need to split variables into real and imaginary parts, resulting in a compact and efficient implementation.

An adaptive initial step size strategy is designed to balance computational efficiency with tracking robustness. The minimum initial step size $\Delta s_{\text{init}}^{\min}$ is chosen as $\max(10^{-3}, 0.1\Delta\lambda_{veering})$, where $\Delta\lambda_{veering}$ is the minimum eigengap in veering regions across solutions at $s=0$. The adaptive initial step size is then
\begin{equation} \label{eq:Delta s_init}
\begin{aligned}
    \Delta s_{\text{init}} = \max \left[\Delta s_{\text{init}}^{\min} \cdot \min\left(\max(1, 2^{\beta-1}), 10)\right), \Delta s_{\text{init}}^{\max}\right],
\end{aligned}
\end{equation}
where $\beta = \Delta \lambda /\Delta \bar{\lambda}$ is the ratio between the actual eigengap $\Delta \lambda$ to the reference gap $\Delta \bar{\lambda}$ (the larger of the 5$\%$ quantile of all eigengaps and twice the minimum eigengap in veering regions $2\Delta\lambda_{veering}$). This formulation yields larger initial steps for well‑separated modes (large $\beta$) and smaller steps for closely spaced modes, ensuring accurate tracking where needed. During continuation, step size is further adapted based on the local curvature: if the inner product of consecutive normalized tangent vectors $\tau = |\dot{\mathbf{y}}_p^H \dot{\mathbf{y}}_{p+1}|$ falls below a threshold $\bar{\tau}=0.99$, the step size is halved; if it remains above, the step size is increased by a factor of 1.1 (up to a prescribed maximum).

Forward tracking continues until the homotopy parameter exceeds the target value, i.e., $s > 1$. At this point, a backward path tracking step is performed with step size $\Delta s = 1 - s$ to return precisely to $s=1$, yielding the exact complex eigensolution $(k_{j,1}, \mathbf{q}_{j,1})$ for the viscoelastic waveguide at the given frequency $\omega_j$. This overshoot-and-refine strategy ensures that the solution at $s=1$ is obtained with full Newton accuracy rather than relying on interpolation from nearby points. Since all homotopy paths are independent, making the method embarrassingly parallel.

The key solutions obtained from Stage 1 consist of a real wavenumber and a frequency that is real in theory but may carry a negligible imaginary component due to numerical round-off. Before homotopy path tracking, any imaginary part of the frequency is discarded, retaining only the physical excitation frequency $\omega = \mathrm{Re}(\omega_{j,0})$. The pair $(\omega, k_{j,0}, \mathbf{q}_{j,0})$ is then used as an initial guess for the homotopy system \cref{eq:G} at $s=0$. Since this guess does not exactly satisfy \cref{eq:G} due to the discarded imaginary part, a single Newton correction step is performed with $\omega$ held fixed. This calibration updates $k$ and $\mathbf{q}$ to satisfy \cref{eq:G} to within a specified tolerance, effectively transferring the small imaginary part from the frequency into the wavenumber while preserving the real-valued frequency as the physical parameter. The calibrated point then serves as the exact starting condition for forward homotopy tracking.

\subsubsection{Theoretical foundations: exceptional‑point topology, branch continuity, and label inheritance}
\label{sec:theory_guarantee}

The observable veering or crossing of dispersion branches along the real frequency axis is a direct projection of the underlying exceptional‑point (EP) topology. Following \cite{keck_unfolding_2003}, two distinct configurations are possible:
\begin{itemize}
    \item \textbf{Type I}: the two EPs lie on opposite sides of the real frequency axis. A frequency scan yields real‑part veering and imaginary‑part crossing; the physical mode identity is preserved without label exchange.
    \item \textbf{Type II}: the two EPs lie on the same side of the real frequency axis. The real parts cross while the imaginary parts veer; physical continuity of the mode shapes then requires exchanging the mode labels at the crossing.
\end{itemize}
This topological classification is the key to deciding whether the transported identities can be inherited automatically or must be corrected.

The proposed homotopy framework links the elastic (Hermitian) limit to the target viscoelastic state through a material attenuation parameter \(s\in[0,1]\). Two theoretical pillars guarantee the continuity of modal identities along this path.
\begin{enumerate}
    \item \textbf{Branch identity continuity.} If the real‑parameter path \(s\in[0,1]\) avoids all EPs in the complex \(s\)-plane, analytic perturbation theory~\cite{kato_perturbation_1995} ensures that every eigenpair can be uniquely continued from \(s=0\) to \(s=1\). This establishes a one‑to‑one \emph{branch identity} between the elastic and viscoelastic solutions, irrespective of the EP topology at the target state.
    \item \textbf{Physical label inheritance.} At the elastic anchor (\(s=0\)), the EPs are symmetric about the real frequency axis, i.e., the system is always in a Type~I configuration. As \(s\) increases, the EPs migrate in the complex frequency plane. If the target state (\(s=1\)) remains Type~I, the physical labels assigned at the anchor are automatically inherited. If the migration causes an EP to cross the real axis, the target becomes Type~II; branch continuity still holds, but the labels must be swapped to reflect the correct physical crossing pattern.
\end{enumerate}

Two practical remarks are in order. First, the parameter space of the homotopy (\(s\)-plane at fixed real frequency) is independent of the observation space (complex frequency plane at \(s=1\)). The migration of EPs in the frequency plane does not imply the presence of an EP on the real interval \(s\in[0,1]\). Second, the predictor–corrector algorithm provides a built‑in safeguard: if the continuation path approaches an EP, the Jacobian becomes ill‑conditioned, causing the corrector to diverge or the step size to fall below a preset minimum, and the algorithm terminates without returning a (potentially incorrect) solution. Successful completion therefore numerically certifies that the path was EP‑free.

Finally, as explicit EP tracking in the complex frequency plane demands dedicated algorithms beyond the scope of this work, the Type II transition cannot be predicted a priori here. Therefore, we introduce two post‑hoc diagnostic signatures in \cref{sec:robustness_violated}—an extremely sharp imaginary‑part crossing and a marked discrepancy between spectral group velocity and energy flux velocity—that reliably indicate where a label swap is necessary. These diagnostics, combined with the topological rule above, give the overall framework both rigor and practical robustness.

\subsubsection{Algorithmic Summary}
\label{sec:algorithm}
\begin{figure}[htb]
\centering
\includegraphics[width=0.9\columnwidth]{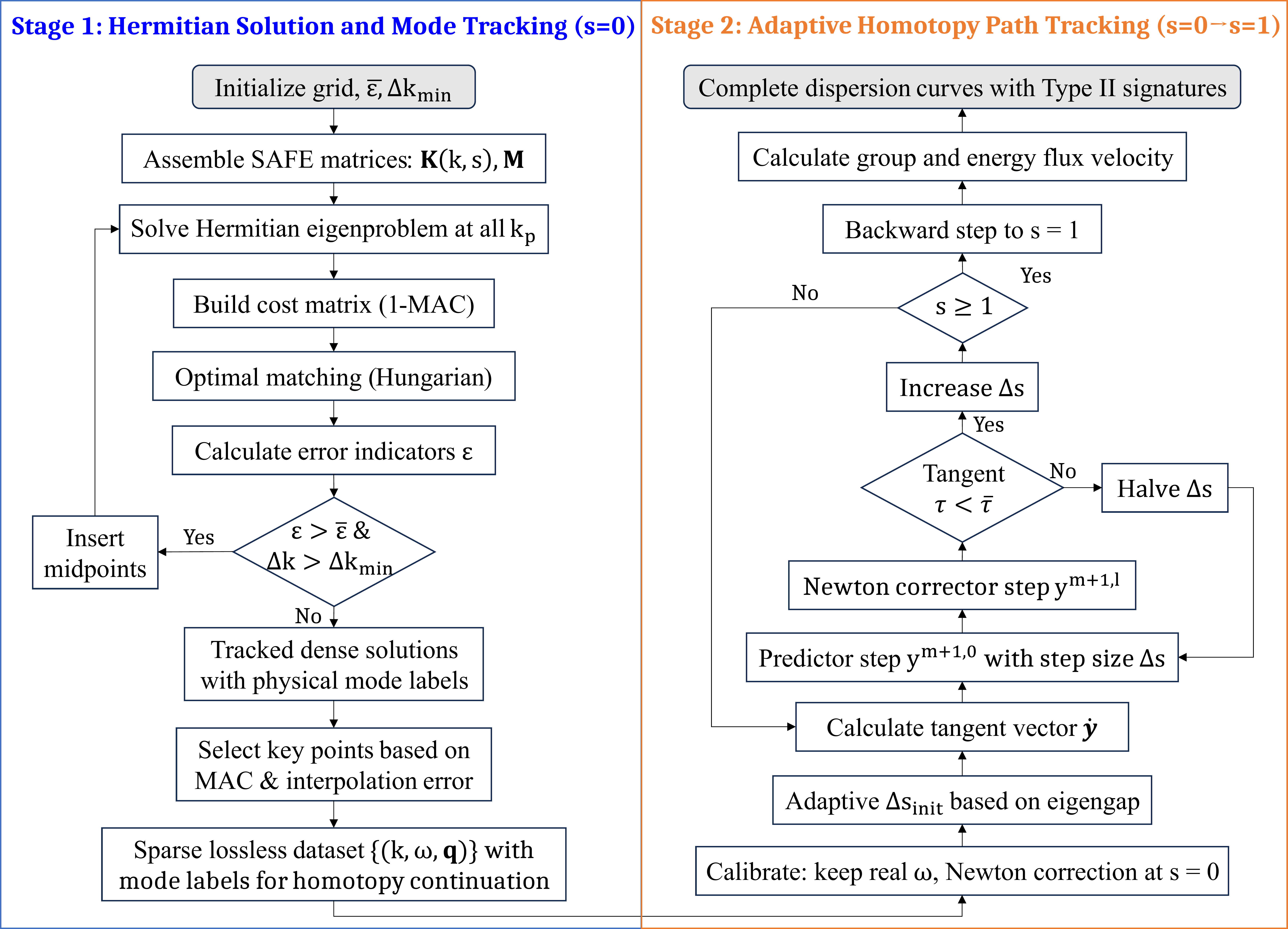}
\caption{Flowchart of the adaptive homotopy continuation framework. The process comprises two stages: (1) Hermitian solution and mode tracking at the lossless stage ($s=0$), and (2) sparse homotopy path tracking ($s=0 \to s=1$) via a predictor-corrector algorithm with adaptive step-size control.}
\label{fig:adap_alg}
\end{figure}

The complete adaptive homotopy continuation framework is illustrated schematically in \cref{fig:adap_alg}. Algorithmically, the framework comprises two principal stages: (i) Hermitian solution and mode tracking at the lossless stage ($s=0$, and (ii) sparse homotopy path tracking from $s=0$ to $s=1$ via a predictor-corrector algorithm with adaptive step-size control. The accurate modal dataset produced by Stage 1 serves as the starting point for the sparse homotopy mapping in Stage 2, ensuring that mode identities are correctly preserved throughout the transition to the viscoelastic regime.

\begin{itemize}
\item \textbf{Stage 1: Hermitian Solution and Mode Tracking ($s=0$)}
    \begin{enumerate}
    \item \textbf{Initialization.} Define the wavenumber range $[k_{\min}, k_{\max}]$, error tolerance $\bar{\varepsilon}$, minimum step size $\Delta k_{\min}$ and the initial grid size $N_0$. 
    
    \item \textbf{Assemble SAFE matrices.} Construct $\mathbf{K}_1$, $\mathbf{K}_2$, $\mathbf{K}_3$, and $\mathbf{M}$ using stiffness tensor $\mathbf{C}(s=0)$.

    \item \textbf{Parallel eigenvalue solution}. At each wavenumber $k_p$, solve the Hermitian eigenproblem $\mathbf{K}(k)\mathbf{q} = \omega^2 \mathbf{M} \mathbf{q}$.
    
    \item \textbf{Optimal mode assignment} For each interval, build a cost matrix based on MAC dissimilarity and apply the Hungarian algorithm to obtain the globally optimal mode matching. Compute the interval error indicator $\varepsilon$. 
    
    \item \textbf{Adaptive refinement.} If $\varepsilon > \bar{\varepsilon}$ and $\Delta k > \Delta k_{\min}$, insert a grid point at the interval midpoint and repeat steps 3–4 until convergence.
    
    \item \textbf{Key point selection.} From the converged dense dataset, select a sparse subset of “key” points satisfying MAC and interpolation error conditions, retaining higher density in veering regions while reducing points where mode shapes evolve slowly.
    \end{enumerate}

\item \textbf{Stage 2: Adaptive homotopy path tracking ($s=0 \to s=1$)}
    \begin{enumerate}
    \item \textbf{Calibration.} For each key point, retain only the real part of the frequency and perform one Newton correction at $s=0$ to obtain a consistent starting point for homotopy.
    \item \textbf{Adaptive initial step size.} Determine $\Delta s_{\text{init}}$ based on the local eigengap: larger steps for well‑separated modes, smaller steps for closely spaced modes (veering regions).
    \item \textbf{Predictor-corrector tracking.} For each calibrated starting point, trace the homotopy path using:
        \begin{itemize}
        \item Tangent prediction with step size $\Delta s$. 
        \item Newton correction to satisfy the extended system.
        \item Adaptive step-size control based on local curvature (inner product of consecutive tangent vectors).
        \end{itemize}
    \item \textbf{Backward refinement.} Once $s>1$ is reached, perform a final backward step to $s=1$ to obtain the exact complex eigensolution.
    \item \textbf{Parallel execution.} All homotopy paths are independent and can be executed in parallel.
    \item \textbf{Group and energy flux velocities computation.} Compute the group and energy flux velocities in parallel, as defined in \cref{eq:group_velocity} and \cref{eq:energy_velocity}, respectively.
    \end{enumerate}
\end{itemize}

The output of Stage 2 is a complete set of complex dispersion curves for the viscoelastic waveguide. Because mode tracking was performed in the Hermitian regime (Stage 1) and the homotopy mapping preserves modal identities, the final curves inherit accurate modal connectivity even in the presence of mode veering, crossing, or degeneracy in the non‑Hermitian system.

\section{Numerical validation}
\label{sec:validation}
This section presents a comprehensive numerical validation of the proposed adaptive homotopy continuation framework. The objective is to demonstrate the algorithm’s capability to accurately and robustly compute dispersion curves with correct mode tracking. 

\subsection{Validation setup and overview of test examples}
\label{sec:validation_setup}
Five numerical examples are selected, covering a diverse range of viscoelastic waveguide geometries and material configurations: two symmetric laminated composite plates, two unsymmetric laminated composite plates, and an L‑shaped aluminium bar. All waveguides are assumed to satisfy traction‑free boundary conditions.

For the laminated composite plates, the performance of the proposed method is benchmarked against the latest version of the freely available Dispersion Calculator (DC) toolbox v3.1 \cite{huber_dispersion_2024, huber_stiffness_2024}. DC computes dispersion curves in viscoelastic layered media using a numerically stable stiffness‑matrix formulation combined with complex wavenumber root‑searching and a predictor–corrector mode‑tracking strategy. Viscoelasticity is incorporated through complex‑valued elastic constants, yielding a complex global stiffness matrix assembled layer‑by‑layer. The dispersion relation is obtained from the characteristic condition $\det \mathbf{K}(k,\omega)=0$, solved via systematic search in the complex $k$-plane followed by local refinement. To construct continuous dispersion curves, DC employs a predictor–corrector continuation in frequency, using previously obtained solutions as initial guesses. Mode continuity is maintained by evaluating proximity in the complex wavenumber space and, when necessary, by comparing mode shapes. These techniques enable DC to robustly capture propagating, attenuated, and leaky guided‑wave modes in viscoelastic plates.

For all numerical examples, the system matrices $\mathbf{K}_j$ ($j=1,2,3$) and $\mathbf{M}$ are normalized using a characteristic length $a$ and a characteristic shear wave velocity $c_T = 3000$ m/s. Consequently, all reported quantities (wavenumber, frequency, phase velocity) are dimensionless. This normalization ensures numerical stability and facilitates consistent parameter selection across examples. The material properties, symmetry, stacking sequences, and characteristic length $a$ for each example are summarized in \cref{tab:params_examples}.
\begin{table}[ht]
\centering
\caption{Material properties, symmetry, stacking sequences, and characteristic length $a$ for the numerical examples.} \label{tab:params_examples}
\small{
\begin{tabular}{llllll}
\toprule
Example & Symmetry & Material & Isotropy & Stacking sequence & Length a \\
\midrule
Laminate Sym1
& Symmetric
& CFRP (Hernando)
& Orthotropic
& $[0, 90, 45, -45]_{2s}$
& \multirow{5}{*}[0ex]{\begin{tabular}{@{}l@{}} Half laminate \\ thickness  \end{tabular}} 
\\
Laminate Sym2
& Symmetric
& CFRP (Castaings)
& Orthotropic
& $[0, 45, -45, 90]_{2s}$ 
& \\
Laminate UnSym1
& Unsymmetric
& CFRP (Hernando)
&  Orthotropic
& $[0, 90, 45, -45]_{4}$
& \\
Laminate UnSym2
& Unsymmetric
& CFRP (Hernando)
& Orthotropic
& \begin{tabular}{@{}l@{}} $[0,15,-15,30,$\\ $-30,45,-45,90]_{2}$ \end{tabular} 
& \\
L-shaped bar 
& Unsymmetric
& Aluminium
& Isotropic 
& -
& Half short leg \\
\bottomrule
\end{tabular}} 
\end{table}

Laminates Sym1 ($[0, 90, 45, -45]_{2s}$) is a symmetric orthotropic laminate (Hernando lamina) that exhibits symmetry‑protected crossings in the elastic state ($s=0$). This example illustrates how these crossings evolve when hysteretic damping is introduced in the viscoelastic state ($s=1$), and provides a direct comparison between the proposed HC method and DC in a configuration where mode veering between two S modes challenges conventional tracking algorithms.

To challenge the framework under more complex modal interactions, we consider two unsymmetric laminates, UnSym1 and UnSym2, both using the Hernando lamina. UnSym1 has a stacking sequence $[0, 45, -45, 90]_{4}$, while UnSym2 adopts the more unbalanced layup $[0,15,-15,30,-30,45,-45,90]_{2}$. The absence of through‑thickness symmetry eliminates symmetry‑protected crossings, leading to pervasive mode veering at $s=0$ that persists at $s=1$.

A second symmetric laminate, Sym2, is constructed using the Castaings lamina with the quasi‑isotropic layup $[0, 45, -45, 90]_{2s}$. With a loss factor of approximately 0.02—substantially larger than that of the Hernando lamina ($\eta \approx 0.003$)—this example tests whether the framework retains its accuracy at elevated damping levels where the exceptional point topology remains Type I, as is typically the case for symmetric laminates with larger veering gaps.

The viscoelastic properties of the Hernando and Castaings laminae are taken from the DC material library \cite{huber_dispersion_2024}, corresponding to the entries \texttt{CarbonEpoxy\_Hernando\_2015\_Viscoelastic} and \texttt{CarbonEpoxy\_Castaings\_Viscoelastic}, respectively. The elastic stiffness matrices of both laminae are listed in \cref{tab:stiffness_matrices}. For anisotropic laminates with component-wise damping, the effective loss factor is defined as the Frobenius-norm ratio $\eta = \| \mathbf{C}'\|_{F}/\| \mathbf{C}''\|_{F}$. For the Hernando lamina this yields $\eta \approx 0.003$ , whereas the Castaings lamina—whose damping is nearly isotropic across all stiffness components—has $\eta = 0.02$.

\begin{table}[ht]
\centering
\caption{Elastic stiffness matrix components (Voigt notation) of the two CFRP laminae used in the composite plate examples. Values are given in GPa. The density of both lamina is 1500 $\mathrm{Kg/m^3}$.}
\label{tab:stiffness_matrices}
\footnotesize{
\begin{tabular}{llllllllll}
\toprule
Material 
& $C_{11}$
& $C_{12}$
& $C_{13}$
& $C_{22}$
& $C_{23}$
& $C_{33}$
& $C_{44}$
& $C_{55}$
& $C_{66}$ \\
\midrule
Hernando 
& 132+0.4i
& 6.9+0.001i
& 5.9+0.016i
& 12.3+0.037i
& 5.5+0.021i
& 12.1+0.043i
& 3.32+0.009i
& 6.21+0.015i
& 6.15+0.02i\\
Castaings 
& 125+2.5i
& 6.3+0.126i
& 5.4+0.108i
& 14+0.28i
& 7.1+0.142i
& 14+0.28i
& 3.45+0.069i
& 5.4+0.108i
& 5.4+0.108i \\
\bottomrule
\end{tabular}}
\end{table}

The L‑shaped aluminium bar is selected to demonstrate the applicability of the proposed framework to general two‑dimensional cross‑sections. 
This example was previously investigated in the purely elastic case in references \cite{maruyama_continuation_2025} and \cite{xiao_rigorous_2026}, which provide details of its geometric dimensions. 
The material properties of aluminium are Young's modulus $E = 70$~GPa, Poisson's ratio $\nu = 0.33$, and density $\rho = 2700$~kg/m$^3$. 
Viscoelastic effects are introduced through complex‑valued Lamé constants:
\begin{equation}
    \begin{aligned}
    \lambda^{\ast} = \lambda\,(1 + i\eta_{\lambda}), \qquad \mu^{\ast} = \mu\,(1 + i\eta_{\mu}),
     \end{aligned}
\end{equation}
with loss factors $\eta_{\lambda}=10^{-4}$ and $\eta_{\mu}=10^{-3}$. 
For polycrystalline aluminium, the loss factor measured by low‑frequency structural resonance techniques is typically of the order of $10^{-4}$ \cite{cremer_structure-borne_1973,granick_material_1965}. 
However, guided wave modes such as Lamb waves involve strong shear components, and higher‑frequency ultrasonic characterisations reveal that the attenuation coefficient of the shear mode is commonly several times larger than that of the longitudinal mode \cite{ono_dynamic_2020}. 
A shear loss factor one order of magnitude above the dilatational loss factor is therefore adopted here, which is also a common assumption in time‑domain viscoelastic modelling of metallic structures.

The adaptive homotopy continuation algorithm is configured with the following parameters, selected based on extensive preliminary experiments:
\begin{itemize}
    \item \textbf{Overall settings:}
        \begin{itemize}
        \item Wave propagation direction: for laminates, dispersion curves are computed for direction $\phi=0$ (material coordinate system of the lamina).
        \item Lamina thickness: each lamina layer has a thickness of 0.25 mm.
        \item Characteristic wave velocity: $c_T = 3000$ m/s for all examples.
        \item Maximum normalized frequency: for laminates, set to a value corresponding to approximately 5000 $\mathrm{kHz \cdot mm}$ (frequency $\cdot$ thickness), covering typical SHM/NDT frequency ranges.
        \item Mesh discretization: for laminates, Gauss–Lobatto–Legendre (GLL) elements with GLL quadrature are used. Each lamina is discretized into two GLL elements of order 5, ensuring adequate through‑thickness resolution for high‑order guided wave modes up to the target frequency–thickness product. For two‑dimensional cross‑sections, nine‑node quadrilateral (quad9) elements with Gaussian quadrature are employed.
        \item Matrix type: all SAFE matrices ($\mathbf{K}_1,\mathbf{K}_2,\mathbf{K}_3,\mathbf{M}$) are stored and processed as sparse matrices to enhance computational efficiency and reduce memory footprint.
        \item Parallel computing: eight concurrent jobs are used for both solving the Hermitian eigenvalue problems (Stage 1) and homotopy continuation (Stage 2).
        \end{itemize}
    \item \textbf{Stage 1:}
        \begin{itemize}
        \item Error tolerance: $\bar{\varepsilon}=0.05$ (the error indicator in \cref{eq:error_indicator} must fall below this value for acceptance).
        \item Minimum wavenumber step: $\Delta k_{\min}=10^{-3}$, preventing infinite refinement due to numerical noise or extremely narrow veering regions.
        \item Initial grid: uniform sampling with 10 points per unit normalized wavenumber ($\Delta k=0.1$), providing a coarse baseline that reveals regions requiring refinement.
        \item Sparse filtering thresholds: MAC threshold $\bar{\zeta}=0.01$, interpolation error threshold $\bar{\gamma}=0.001$, ensuring that points with rapid eigenvector variation are retained.
        \end{itemize}
    \item \textbf{Stage 2:}
        \begin{itemize}
        \item Maximum initial step size: $\Delta s_{\text{init}}^{\max}=0.01$, enabling fast homotopy tracking for key solutions with large eigengaps.
        \item Tangent vector similarity threshold: $\bar{\tau}=0.99$, providing an additional safety margin for homotopy path continuation.
        \end{itemize}
\end{itemize}

In Stage1, at each refinement iteration the error indicator is recomputed for all marked intervals, and the process continues until convergence or until the minimum step size is reached. In Stage 2, the step size $\Delta s$ is adaptively adjusted during predictor–corrector tracking for each solution point. The following subsections present the dispersion curve results for each numerical example, comparing the proposed method against DC. 

\subsection{Symmetric laminate: symmetry-protected crossings}
\label{sec:example_Sym1}
\begin{figure}[htb]
\centering
\includegraphics[width=1.02\columnwidth]{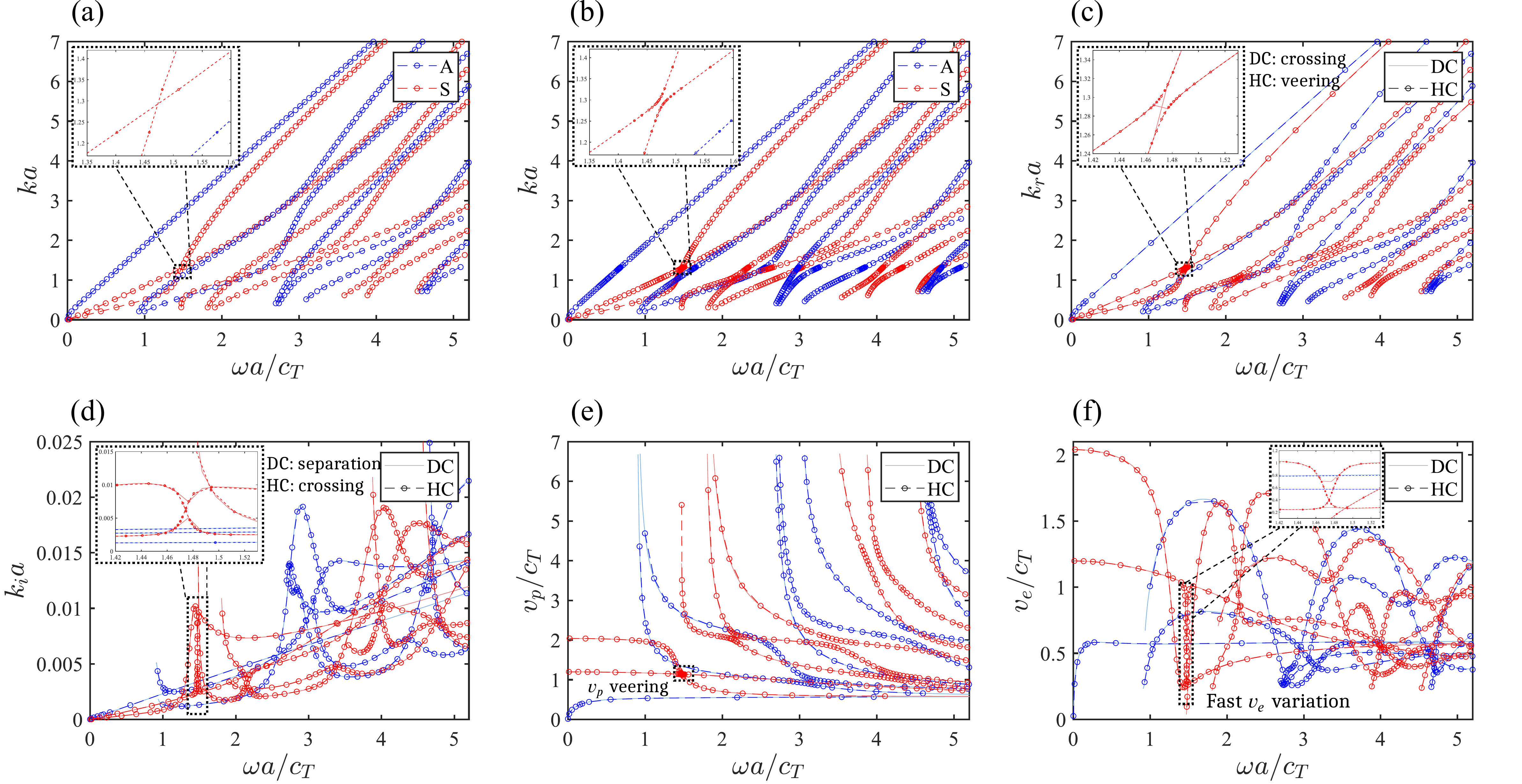}
\caption{Dispersion curves for symmetric laminate Sym1 (A modes: blue, S modes: red). (a)–(b) Elastic state ($s=0$): uniform sampling (a) misidentifies veering as crossing, while adaptive refinement (b) correctly resolves the veering. (c)–(f) Viscoelastic state ($s=1$): comparison between HC (deep dashed) and DC (light solid). (c) Real wavenumber: HC and DC agree except for misidentified veering (black squares). (d) Imaginary wavenumber: DC shows artificial separation, HC correctly captures imaginary-part crossing. (e)–(f) Phase and energy flux velocities: veering appears as a crossing in energy flux velocity (black square).}
\label{fig:hdamp_exa1_sym_CEcast_DCvsHC}
\end{figure}

We first consider a symmetric laminate to demonstrate the proposed framework on a problem where symmetry-protected crossings arise naturally. For validation, we benchmark our method against the freely available Dispersion Calculator (DC) \cite{huber_dispersion_2024, huber_stiffness_2024}, which determines dispersion curves in viscoelastic layered media using a stiffness‑matrix formulation combined with complex wavenumber root‑searching and a predictor–corrector mode‑tracking strategy. \Cref{fig:hdamp_exa1_sym_CEcast_DCvsHC} presents the dispersion curves computed using both the proposed homotopy continuation (HC) method and DC for the symmetric laminate Sym1 (in \cref{tab:params_examples}) with stacking sequence $[0, 90, 45, -45]_{2s}$, consisting of Hernando lamina (in \cref{tab:stiffness_matrices}).

\textbf{Elastic state ($s=0$).} \Cref{fig:hdamp_exa1_sym_CEcast_DCvsHC}(a) and (b) show the frequency–wavenumber curves at the elastic limit. The symmetric layup gives rise to widely observed symmetry-protected crossings between antisymmetric (A) and symmetric (S) modes, which are preserved regardless of the sampling strategy. \Cref{fig:hdamp_exa1_sym_CEcast_DCvsHC}(a) employs uniform wavenumber sampling; in the veering region between two S modes, the eigenvectors change rapidly over a narrow interval. With uniform sampling, this rapid variation is insufficiently resolved, leading to inaccurate mode tracking where the veering is misidentified as a crossing. The adaptive refinement strategy shown in \cref{fig:hdamp_exa1_sym_CEcast_DCvsHC}(b) automatically adds sampling points where the eigenvector variation is high, thereby correctly resolving the veering and restoring accurate modal connectivity. Further theoretical justification and numerical validation of this adaptive strategy are provided in our previous work \cite{xiao_rigorous_2026}.

\textbf{Viscoelastic state ($s=1$).} \Cref{fig:hdamp_exa1_sym_CEcast_DCvsHC}(c)–(f) compare the HC and DC results at the target viscoelastic state. The results of the two methods are nearly identical except in regions highlighted by black squares, which correspond to a veering region between two S modes. In this region, DC misidentifies the veering as a crossing, leading to a persistent mode tracking error. The DC output in this region exhibits two distinct artifacts indicative of tracking failure. First, the real-wavenumber curves are not smooth: the discrete points show rapid, discontinuous changes in curvature, and the curve segments appear jagged rather than following a smooth veering profile. This lack of smoothness indicates numerical instability in the DC mode‑tracking algorithm, where the solver occasionally jumps between nearby branches and then back, producing kinks and non‑monotonic segments. Second, in the imaginary wavenumber plot (\cref{fig:hdamp_exa1_sym_CEcast_DCvsHC}(d)), the inaccurate DC tracking manifests as an artificial separation between the two branches in $k_i$, whereas the correct HC results exhibit a crossing—the physically expected Type I behavior for hysteretic damping, where real parts veer and imaginary parts cross.

Crucially, the HC results show that the real‑part veering gap remains moderately large (the two S modes never approach extremely close), and the imaginary‑part crossing occurs over a broad frequency interval rather than being compressed into a very narrow band. Such a configuration is characteristic of a system that is still well within the Type I regime, where the two EPs remain far from the real axis. Consequently, the erroneous ``crossing'' reported by DC is not a sign of an impending topological transition; it is purely a tracking artifact. The fact that DC misinterprets a clear Type I veering as a crossing, even in this mildly damped and well‑separated case, underscores the fragility of its post‑processing heuristics.

Comparing \cref{fig:hdamp_exa1_sym_CEcast_DCvsHC}(b) and (c), the introduction of hysteretic damping does not significantly alter the real-part dispersion curves; the veering and symmetry-protected crossings remain present, though slight frequency shifts are observed. This behavior is expected because the frequency-independent damping preserves the analytic structure of the eigenvalue problem along the frequency axis.

\Cref{fig:hdamp_exa1_sym_CEcast_DCvsHC}(e) and (f) show the corresponding phase velocity and energy flux velocity curves. In the phase velocity plot, the veering appears as a characteristic approach and recession of the two curves, maintaining a positive gap. In the energy flux velocity plot, however, the same veering manifests as an apparent crossing: the two modes rapidly exchange their order and positions, creating the illusion of an intersection. This energy flux velocity behavior, highlighted by the black square, further illustrates the subtlety of mode tracking in veering regions and underscores the necessity of a method that reliably preserves modal identity.

\begin{figure}[htb]
\centering
\includegraphics[width=0.8\columnwidth]{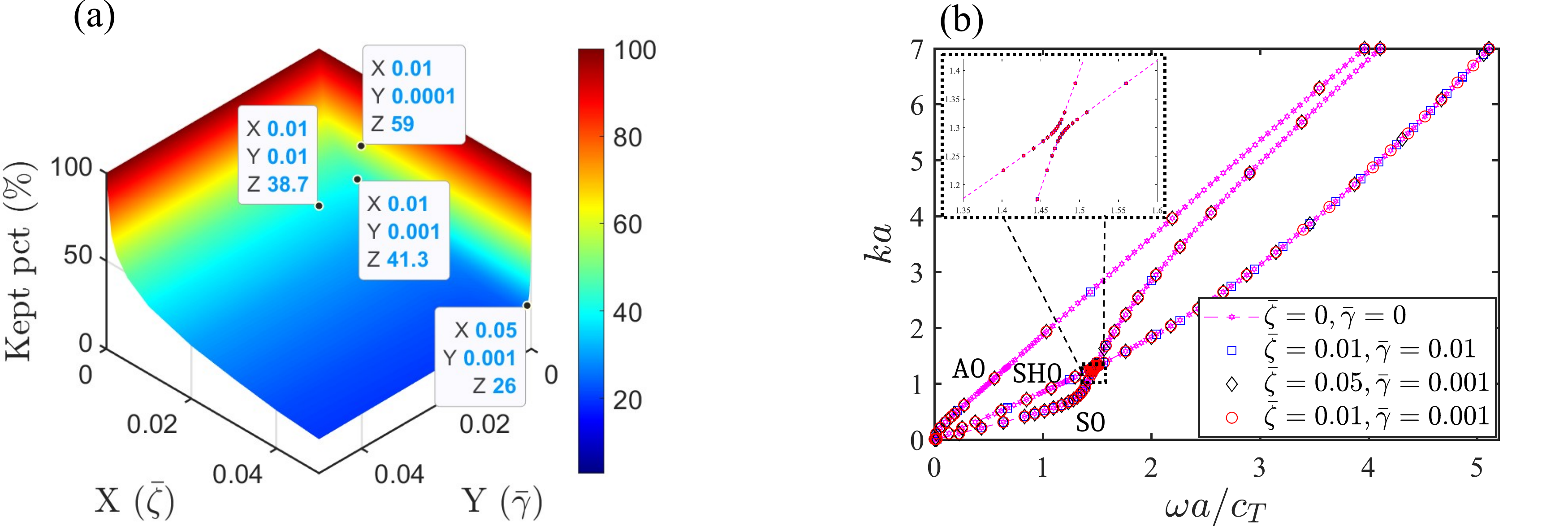}
\caption{Selection of key elastic (lossless) solutions ($s=0$) for homotopy continuation. (a) Surface plot of the percentage of retained solutions as functions of the filtering parameters $\bar{\zeta}$ (MAC parameter) and $\bar{\gamma}$ (interpolation error parameter). (b) Solution points for modes A0, SH0, and S0 before and after filtering for three parameter combinations; all candidate points are indicated as ($\bar{\gamma}=0, \bar{\zeta}=0$) in the legend}
\label{fig:hdamp_exa1_sym_solutionfilter}
\end{figure}

To illustrate the selection of key elastic solutions at the lossless state ($s=0$) for homotopy continuation, \cref{fig:hdamp_exa1_sym_solutionfilter} demonstrates how the two filtering parameters—the MAC parameter $\bar{\zeta}$ and the interpolation error parameter $\bar{\gamma}$ —control the solution filtering process.

\Cref{fig:hdamp_exa1_sym_solutionfilter}(a) presents a surface plot of the percentage of retained solutions as a function of $\bar{\zeta}$ and $\bar{\gamma}$. As expected, the retained percentage decreases with increasing values of either parameter. Notably, the reduction rate is steeper along the interpolation error parameter $\bar{\gamma}$ in the range [0,0.001] than along the MAC parameter $\bar{\zeta}$. Fixing $\bar{\zeta}= 10^{-2}$, increasing $\bar{\gamma}$ from 0 to $10^{-4}$ and $10^{-3}$ yields retained percentages of 59.0$\%$ and 41.3$\%$, respectively. Further increase of $\bar{\gamma}$ to $\bar{\gamma}= 10^{-2}$ does not substantially decrease the retained percentage, as $\bar{\zeta}= 10^{-2}$ already guarantees a minimum retention of 38.1$\%$ for this symmetric laminate over the considered frequency range. Conversely, fixing $\bar{\gamma}=10^{-3}$ guarantees a minimum retained percentage of 20.8$\%$ regardless of $\bar{\zeta}$.

\cref{fig:hdamp_exa1_sym_solutionfilter}(b) further illustrates the distribution of solution points for three modes (A0, SH0, S0) before and after filtering using three representative parameter combinations. All candidate points (before filtering) are labeled in the legend as $\bar{\zeta}=0, \bar{\gamma}=0$. For all three combinations—($\bar{\zeta}=0.01, \bar{\gamma}=0.01$), ($\bar{\zeta}=0.05, \bar{\gamma}=0.001$), and ($\bar{\zeta}=0.01, \bar{\gamma}=0.001$)—the points located in the mode veering regions are consistently retained. This behavior is by design: both parameters preferentially retain solutions where eigenvalues or eigenvectors undergo significant variation. Either $\bar{\zeta}=0.01$ or $\bar{\gamma}=0.001$ alone proves sufficient to preserve all veering points. In contrast, the majority of discarded points lie in flat curve regions, where the eigenvalue and eigenvector vary slowly. For instance, in the flat portion of the A0 mode beyond a normalized frequency of 0.5, the three parameter combinations retain only 2, 4, 4 point(s), respectively. The overall retention percentages for all modes are 38.6 38.6$\%$, 26.0$\%$, 41.3$\%$ for the three combinations. To achieve a well‑balanced representation of the final homotopy solutions at $s=1$, the last combination ($\bar{\zeta}=0.01, \bar{\gamma}=0.001$) is selected for this study. 

\begin{figure}[htb]
\centering
\includegraphics[width=0.8\columnwidth]{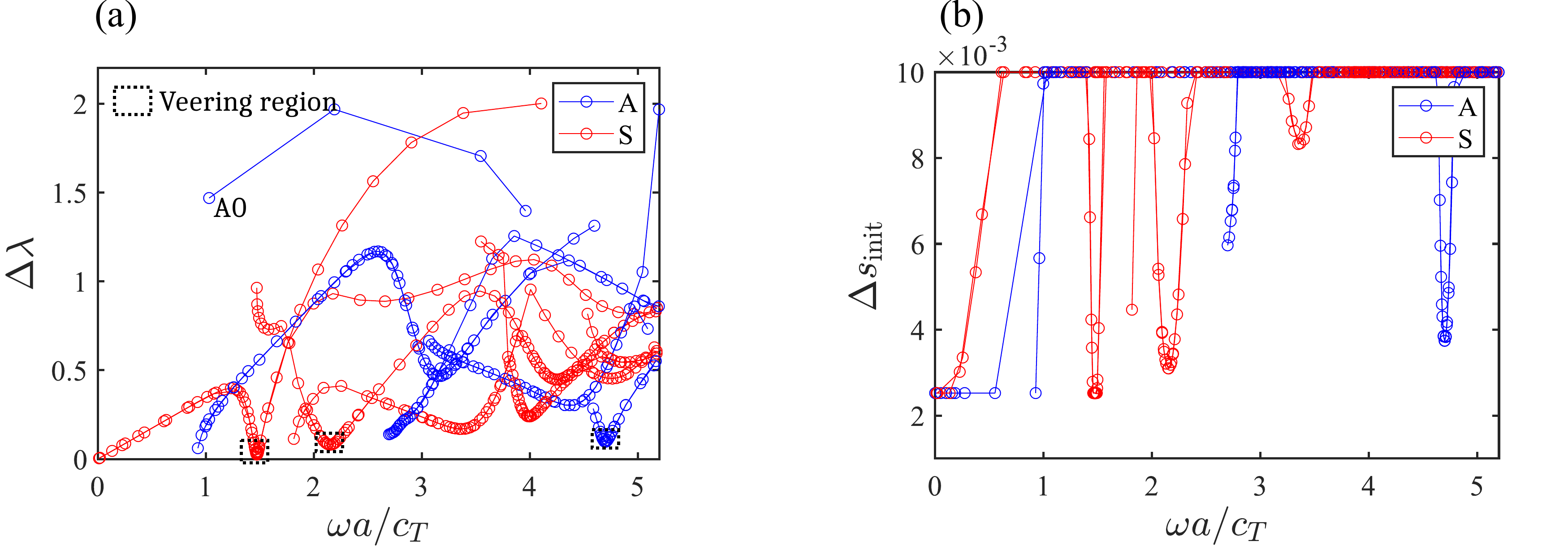}
\caption{Eigengap of selected key solutions and adaptive determination of the initial step size $\Delta s_{\text{init}}$  for homotopy tracking from $s=0 \to s=1$. (a) Eigengap distribution indicating the risk of mode jumping; (b) Adaptive initial step size based on eigengap, with smaller steps assigned to regions of small eigengap.}
\label{fig:hdamp_exa1_sym_eigengap_ds}
\end{figure}

\cref{fig:hdamp_exa1_sym_eigengap_ds} illustrates the eigengap of the selected key solutions and the corresponding adaptive determination of the initial step size $\Delta s_{\text{init}}$ for each solution prior to homotopy tracking from $s=0 \to s=1$. The initial step size is adaptively determined based on the eigengap within the same wave family, estimated from the refined solutions at $s=0$, , as shown in \cref{fig:hdamp_exa1_sym_eigengap_ds}(a). Local minima of these eigengap curves indicate mode veering regions, where a smaller $\Delta s_{\text{init}}$ is required to ensure accurate tracking. The minimum eigengap across all mode veering, i.e. the minimum of the veering eigengaps $\Delta \lambda_{veering}$, is 0.0252, and the veering threshold is adaptively set to 5$\%$ quantile of all eigengaps, i.e., $\Delta \bar{\lambda} = 0.0622$. Additionally, near cutoff frequencies where a mode lacks a counterpart within the same wave family (e.g., the A0 mode at normalized frequencies below 0.9 in this example), the solutions are also assigned a small $\Delta s_{\text{init}}$.

The adaptive determination of $\Delta s_{\text{init}}$ using \cref{eq:Delta s_init} is presented in \cref{fig:hdamp_exa1_sym_eigengap_ds}(b). For solution points whose eigengap $\Delta \lambda$ does not exceed $\Delta \bar{\lambda}$—including points near cutoff frequencies without an explicit $\Delta \lambda$ —the initial step size is adaptively set to $0.1\Delta \lambda_{veering} = 0.00252 > \Delta s_{\text{init}}^{\min} = 0.001$. As $\Delta \lambda$ increases beyond the threshold, $\Delta s_{\text{init}}$ increases accordingly, up to the prescribed maximum $\Delta s_{\text{init}}^{\max}=0.01$ used in this study.

These results for the symmetric laminate Sym1 demonstrate that the key solution filtering strategy successfully retains points in regions of rapid eigenvector variation, while the adaptively assigned initial step sizes ensure reliable homotopy tracking from $s=0$ to $s=1$. In the viscoelastic state, the veering is correctly resolved by HC as a Type I crossing, whereas DC misidentifies it as a real-part crossing. This contrast highlights the fragility of direct non-Hermitian mode tracking and motivates the more demanding unsymmetric validation cases that follow.


\subsection{Unsymmetric laminate: pervasive Type I behavior}
\label{sec:example_UnSym1}
\begin{figure}[htb]
\centering
\includegraphics[width=0.8\columnwidth]{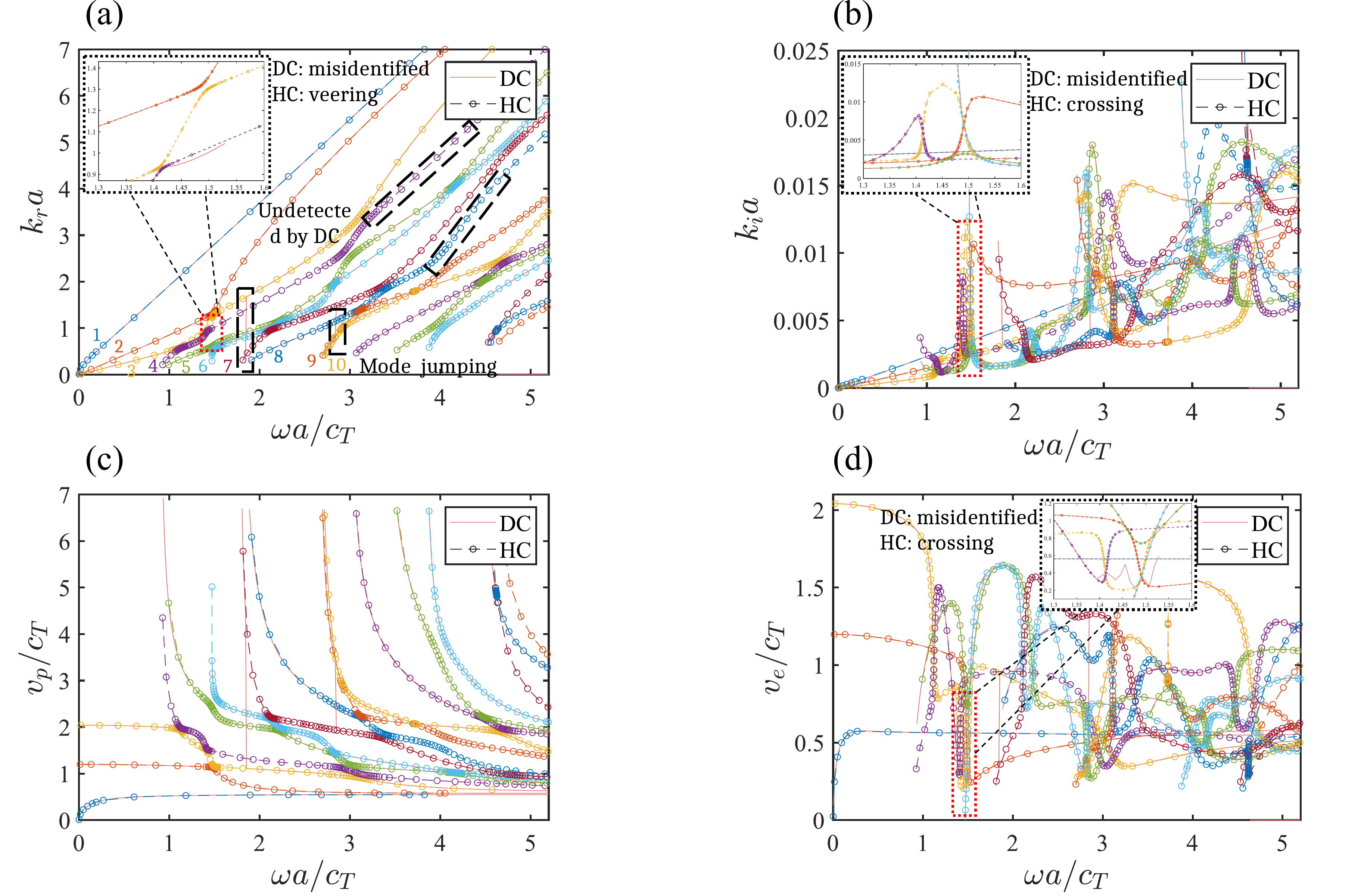}
\caption{Dispersion curves for unsymmetric laminate UnSym1 at the viscoelastic state ($s=1$). Comparison between HC (dashed lines with markers) and DC (light red solid lines). (a) Real wavenumber; (b) Imaginary wavenumber; (c) Phase velocity; (d) energy flux velocity. Red squares in (a) indicate type I crossing regions where DC incorrectly tracks mode crossing.}
\label{fig:hdamp_Unsym_exa1_CEcast_DCvsHC}
\end{figure}

The symmetric laminate examples demonstrate that mode veering poses significant challenges to mode tracking, particularly at the viscoelastic state ($s = 1$). To further assess the robustness of the proposed HC method in handling pervasive veering, we consider an unsymmetric laminate UnSym1 with stacking sequence $[0, 90, 45, -45]_{4}$ (see \cref{tab:params_examples}). \Cref{fig:hdamp_Unsym_exa1_CEcast_DCvsHC} presents the dispersion curves at the viscoelastic state ($s=1$), comparing the HC method with the reference DC. To facilitate direct comparison, each subfigure overlays results from both methods: DC is shown as light red solid lines, while HC is represented by dashed lines with round markers indicating the solution points.

In unsymmetric laminates, the absence of through-thickness symmetry eliminates symmetry-protected crossings, leading to pervasive mode veering in the elastic state ($s=0$) \cite{xiao_rigorous_2026}. This characteristic persists at the viscoelastic state ($s=1$), as shown in the frequency–real wavenumber dispersion plot of \cref{fig:hdamp_Unsym_exa1_CEcast_DCvsHC}(a). The numerical results obtained with the HC method fully corroborate the theoretical analysis presented in \cref{sec:theory_guarantee}: when no topological transition occurs, Type I behavior dominates. Consequently, the veering events observed at $s=0$ manifest as Type I behavior at $s=1$—real‑part veering (\cref{fig:hdamp_Unsym_exa1_CEcast_DCvsHC}(a)) accompanied by imaginary‑part crossing (\cref{fig:hdamp_Unsym_exa1_CEcast_DCvsHC}(b)). Notably, the imaginary‑part crossing in the HC results occurs over a relatively wide frequency interval, indicating that the two EPs remain well separated from the real axis and the system is securely in the Type I configuration. The modal identities are therefore preserved along the frequency axis.

In contrast, the reference DC method exhibits significant difficulties in this pervasive veering environment. As shown in \cref{fig:hdamp_Unsym_exa1_CEcast_DCvsHC}(a), one mode (numbered 4) is entirely undetected. Moreover, the missing modes exacerbate tracking errors: for instance, near the leftmost red square at a normalized frequency of approximately 1.4, mode 3 erroneously jumps to mode 4 in the tracking process, resulting in early termination of that branch. Similarly, mode 7 jumps to mode 3 near frequency 1.9, and mode 9 jumps to mode 7 near frequency 2.8. These spurious mode jumps are indicative of incorrect branch assignments.

The failure of DC stems from its mode‑by‑mode tracking strategy. When the predictor–corrector continuation encounters a veering region where eigenvalues are closely spaced and eigenvectors change rapidly, the algorithm may lose the correct branch. To prevent indefinite stalling, DC employs an internal termination criterion that aborts the tracking of the current mode and proceeds to the next one. This heuristic explains the missing modes (e.g., mode 4) and the observed jumps: once a branch is lost, the subsequent tracking may inadvertently continue on a different branch, leading to mislabeling. Consequently, the DC output contains numerous gaps and artifacts that do not correspond to any physically admissible spectral structure, underscoring the need for the proposed homotopy framework which performs mode tracking once and for all in the Hermitian regime.

The proposed HC method circumvents this issue entirely by performing mode tracking exclusively in the Hermitian regime ($s=0$), where eigenvectors are well‑behaved and the MAC is reliable; correct modal identities are then propagated to $s=1$ along the homotopy path. The theoretical analysis in \cref{sec:theory_guarantee} establishes that, for hysteretic damping and within the regime considered, the branch cut topology is unchanged and the EPs remain in the Type I configuration. Consequently, the identities established at $s=0$ remain valid along the entire real frequency axis at $s=1$. The HC results therefore faithfully represent the underlying physics, free from the spurious mode exchanges and unphysical crossing patterns that plague the DC output.

\subsection{Unsymmetric laminate with more unbalanced stacking sequences}
\label{sec:example_UnSym2}
To further challenge the robustness of the proposed HC method, we construct a more complex unsymmetric laminate, UnSym2, by introducing multiple off-axis plies ($\pm 15^{\circ}, \pm 30^{\circ}$) that intentionally reduce eigengaps and induce dense mode veering, thereby increasing the difficulty of mode tracking. The stacking sequence of UnSym2 is $[0,15,-15,30,-30,45,-45,90]_{2}$.

\begin{figure}[htb]
\centering
\includegraphics[width=0.8\columnwidth]{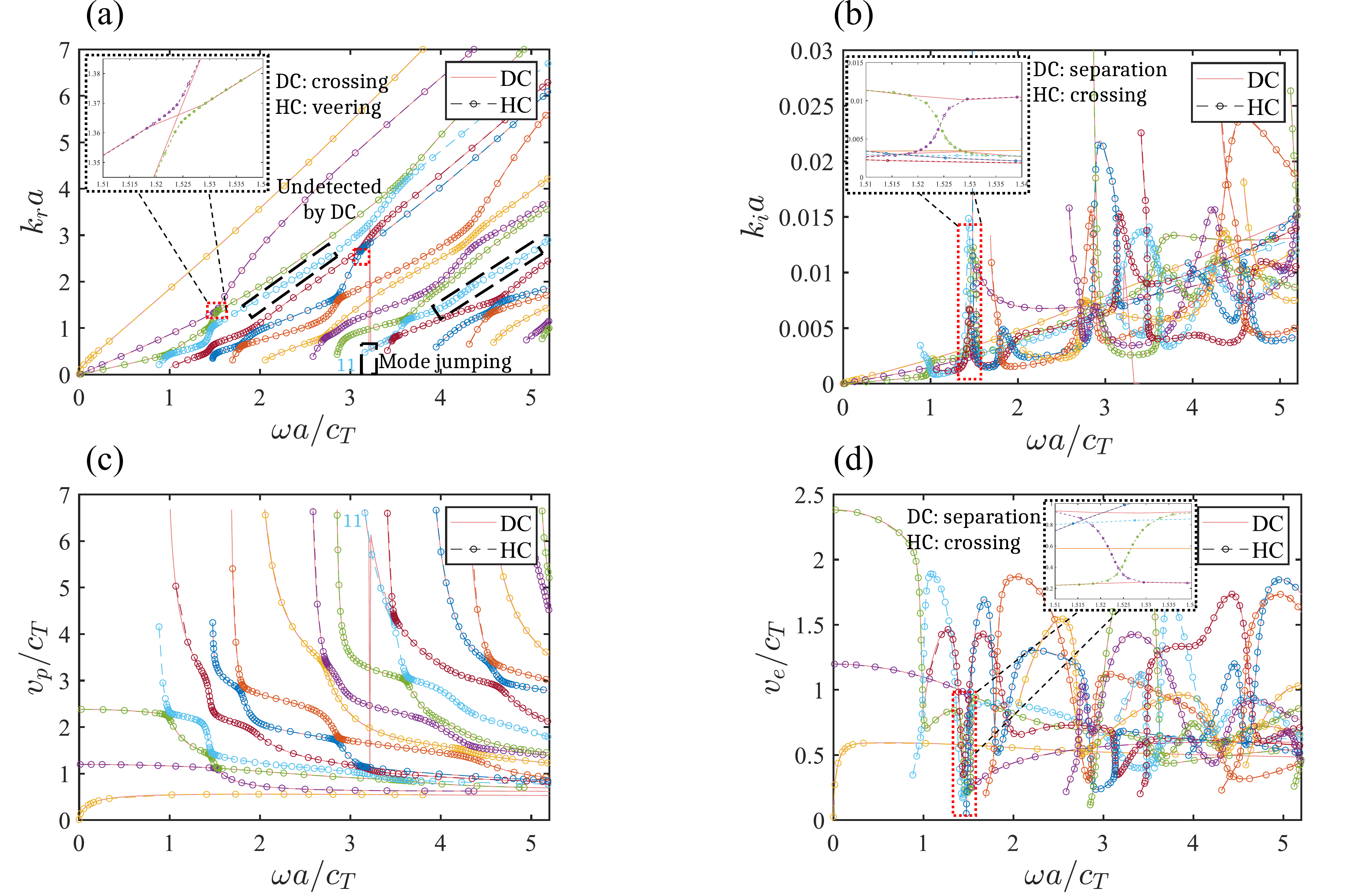}
\caption{Dispersion curves for unsymmetric laminate UnSym2 at the viscoelastic state ($s=1$). Comparison between HC (dashed lines with markers) and DC (light red solid lines). (a) Real wavenumber; (b) Imaginary wavenumber; (c) Phase velocity; (d) energy flux velocity. Red squares in (a) indicate veering regions where DC incorrectly tracks mode crossing.}
\label{fig:hdamp_Unsym_exa2_CEhern_DCvsHC}
\end{figure}

\Cref{fig:hdamp_Unsym_exa2_CEhern_DCvsHC} compares the dispersion curves obtained with the HC method and the reference DC for this laminate at the viscoelastic state ($s=1$). The two methods agree well in regions where modal interactions are weak, but significant discrepancies appear in veering zones, where the real wavenumber curves of DC exhibit artificial crossings (marked by red squares).

The HC results consistently exhibit Type I behavior: the real parts veer (maintaining a positive gap), while the imaginary parts cross over a relatively broad frequency interval. This pattern is physically expected for the low‑damping Hernando lamina ($\eta\approx0.003$). In contrast, DC misinterprets the veering as real‑part crossing and produces imaginary‑part curves that are nearly separated, resembling a gap rather than a veering. For instance, the interaction between modes 2 and 3 is captured by DC, but the imaginary parts appear as two parallel lines with a nearly constant separation—a pattern that, according to Keck crossing theory \cite{keck_unfolding_2003}, would require the two exceptional points to lie on the same side of the real axis and far away from it. Such a configuration is highly implausible for this weakly damped laminate, where the loss factor is only 0.003. The observed artificial separation is therefore a numerical artifact, not a genuine topological feature.

Consistent with previous examples, DC fails to capture all modes, missing one mode entirely. Moreover, for mode 11, DC loses track after a significant jump to another branch; the tracking process terminates prematurely at a normalized frequency of 3.43, yielding only 30 solution points. This premature termination is particularly evident in the phase velocity plot (\cref{fig:hdamp_Unsym_exa2_CEhern_DCvsHC}(c)). These failures underscore the sensitivity of DC to early mode detection and eigengap size: missing modes can trigger spurious jumps, while small eigengaps can lead to erroneous real‑wavenumber crossings.

The HC framework avoids the failures observed in UnSym2 for the same structural reasons detailed in \cref{sec:example_UnSym1}: mode tracking performed once in the Hermitian regime eliminates the detection-omission and branch-jump pathways that afflict DC's sequential strategy. The additional complexity of UnSym2—multiple off-axis plies reducing eigengaps—merely amplifies these failures without introducing new error mechanisms, confirming that the HC advantage is systematic rather than case-specific.

\subsection{Symmetric laminate with elevated damping}
\label{sec:example_Sym2}
\begin{figure}[htb]
\centering
\includegraphics[width=0.8\columnwidth]{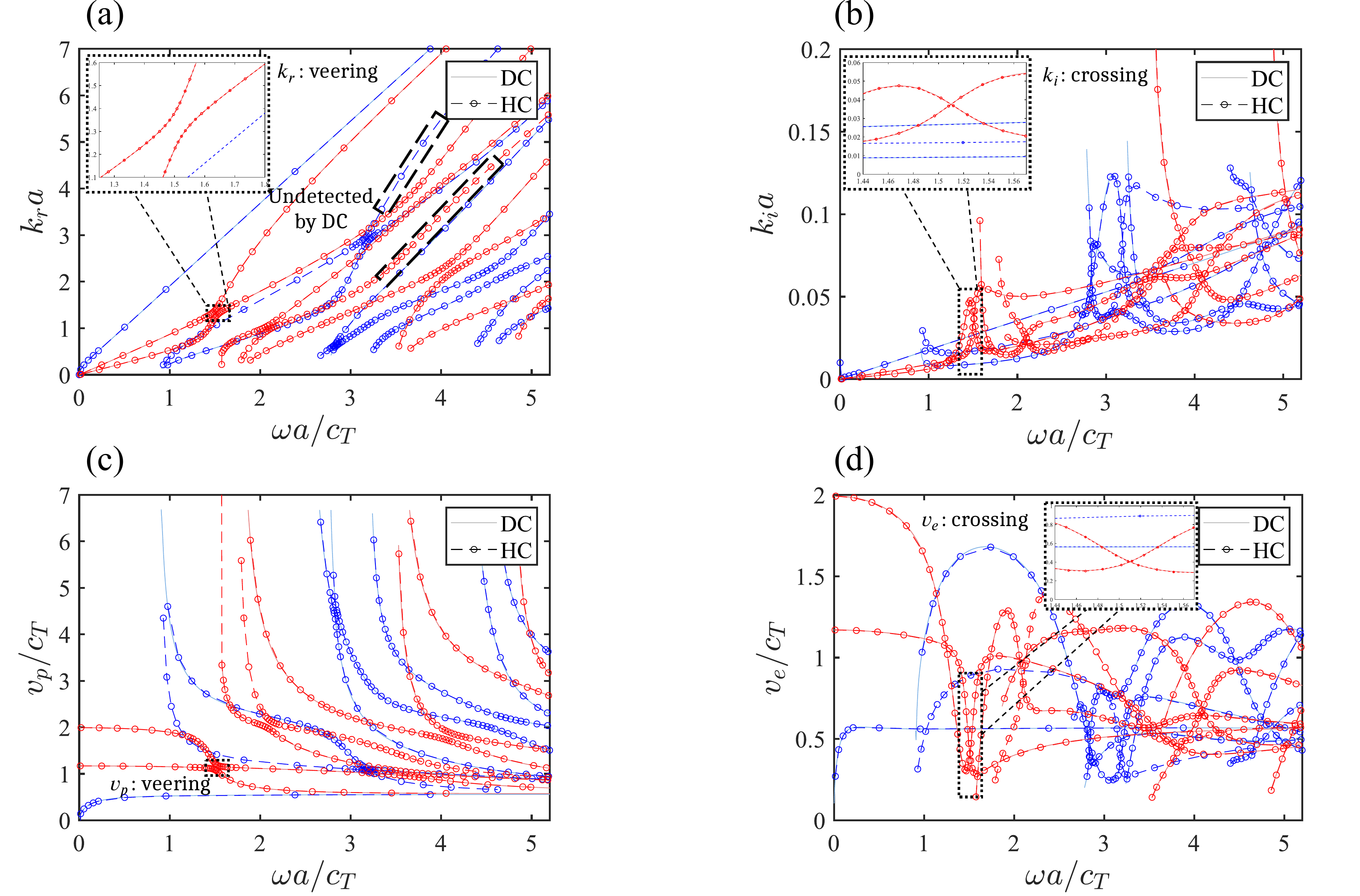}
\caption{Dispersion curves for symmetric laminate Sym2 (A modes: blue, S modes: red) at the viscoelastic state ($s=1$). Comparison between HC (deep dashed) and DC (light solid). (a) Real wavenumber; (b) Imaginary wavenumber; (c) Phase velocity; (d) energy flux velocity. The black squares indicate regions where DC fails.}
\label{fig:hdamp_Sym_exa2_DCvsHC}
\end{figure}

To further evaluate the performance of the proposed HC method under elevated damping, we consider a different lamina material, Castaings (see \cref{tab:stiffness_matrices}), which has a loss factor of approximately 0.02—substantially larger than that of Hernando. This numerical example, denoted Sym2, consists of a symmetric layup $[0, 45, -45, 90]_{2s}$ using the Castaings lamina. \Cref{fig:hdamp_Sym_exa2_DCvsHC} presents the dispersion curves at the viscoelastic state $s=1$, comparing HC with the reference DC results.

Compared to Sym1, Sym2 possesses a relatively large eigengap at the elastic state: the minimum eigengap $\Delta \lambda$ is 0.0867, approximately three times larger than that of Sym1 (0.0252). As derived in our previous work \cite{xiao_rigorous_2026}, a larger eigengap implies a slower eigenvector variation, which facilitates reliable mode tracking even with simpler continuation strategies. Consequently, DC successfully identifies the veering between the two S modes in Sym2, as shown in the zoomed view of \cref{fig:hdamp_Sym_exa2_DCvsHC}(a). Both HC and DC correctly capture the Type I behavior in this veering region: the real parts exhibit a clear veering with a positive gap, the imaginary parts cross over a relatively broad frequency interval, and the energy flux velocity shows the characteristic rapid exchange of mode order (see \cref{fig:hdamp_Sym_exa2_DCvsHC}(b)–(d)). 

Despite this improvement, DC still fails to capture two modes elsewhere in the spectrum, as indicated by the black squares in \cref{fig:hdamp_Sym_exa2_DCvsHC}(a). Even when the detection frequency resolution is reduced to 0.025 kHz, these modes remain undetected. The missing modes and the resulting incomplete dispersion curves underscore the fragility of DC’s post‑processing heuristics, which rely on frequency‑axis continuation and mode shape comparison in the non‑Hermitian domain.

In contrast, the HC framework, by performing mode tracking once in the Hermitian regime ($s=0$) and then propagating identities via homotopy continuation, captures the full spectrum without any missing branches. The larger eigengap in Sym2 also benefits HC by allowing larger initial step sizes and fewer adaptive refinements, further improving computational efficiency. This example demonstrates that the proposed method remains accurate and robust at an elevated damping level of $\eta = 0.02$ for a symmetric laminate, as long as the system stays within the Type I regime—which is confirmed by the broad imaginary crossing and moderate energy flux velocity exchange. The fact that DC, despite correctly handling the veering region, still fails to detect complete modes highlights the fundamental advantage of the HC framework: decoupling mode tracking from the non‑Hermitian solution process ensures that once the elastic modes are correctly identified, the viscoelastic modes are inherited reliably, regardless of damping magnitude.

\subsection{L-shaped bar: arbitrary cross-section}
\label{sec:example_Lbar}
A key advantage of the SAFE formulation over matrix methods such as the SMM lies in its ability to handle arbitrary two-dimensional cross-sections. To demonstrate the applicability of the proposed adaptive HC framework to such general geometries, we consider the L-shaped aluminium bar introduced at the beginning of \cref{sec:validation}. The material properties and geometric parameters are as described therein.

\begin{figure}[htb]
\centering
\includegraphics[width=1.02\columnwidth]{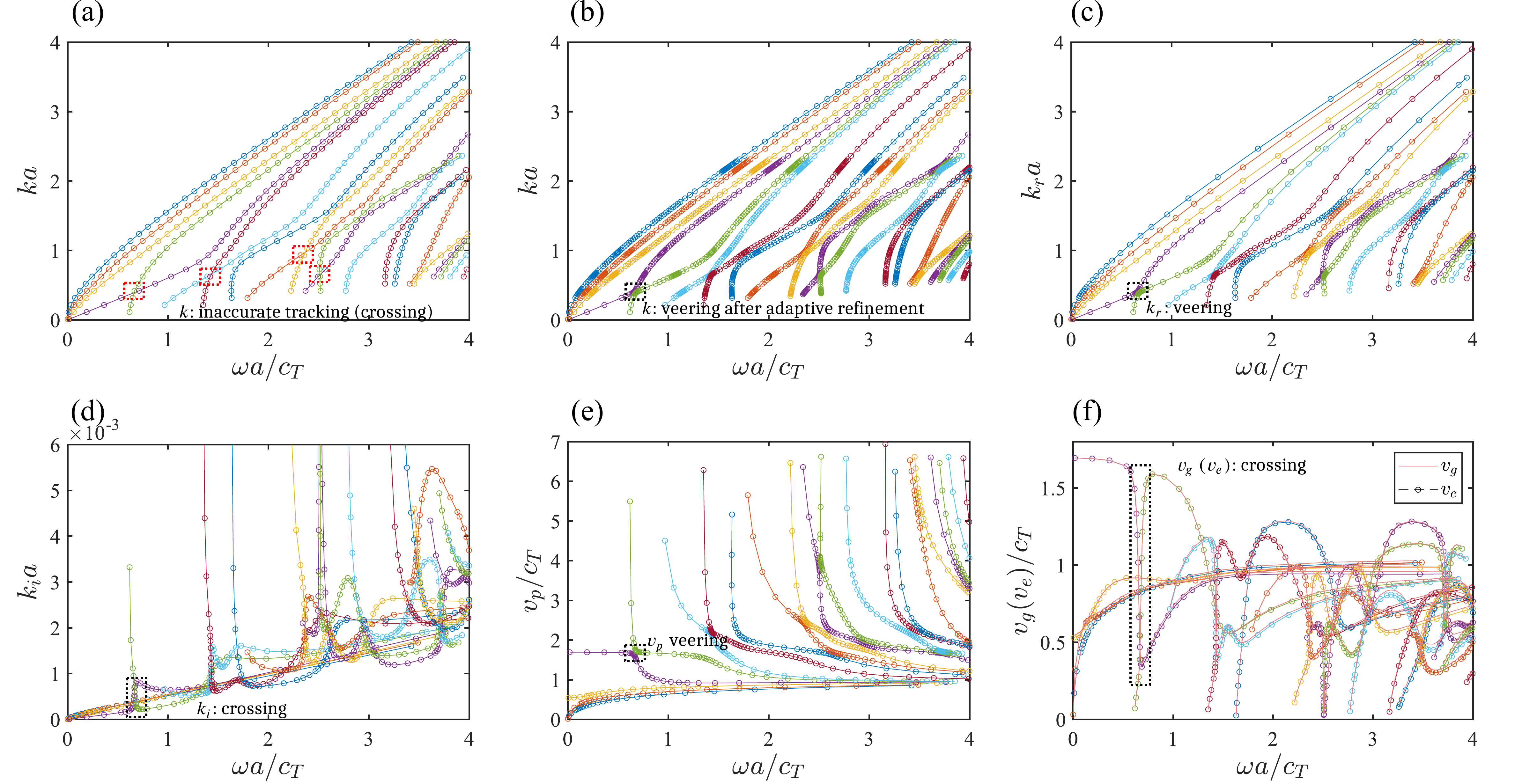}
\caption{Dispersion curves for the L-shaped aluminium bar at the elastic ($s=0$) and viscoelastic ($s=1$) states. (a) Elastic state, uniform sampling ($\Delta k = 0.1$): an erroneous crossing appears (red squares). (b) Elastic state, adaptive refinement: the veering is correctly resolved. (c)–(f) Viscoelastic state: (c) real wavenumber, exhibiting clear veering; (d) imaginary wavenumber, with crossings occurring over a relatively broad frequency interval; (e) phase velocity; (f) spectral group velocity $v_g$ and energy flux velocity $v_e$. A slight offset between the two velocities is visible, with $v_g$ consistently slightly larger than $v_e$; nevertheless, they retain identical structural shapes in the mode-interaction regions, indicating that the target viscoelastic system remains well within the Type I regime.}
\label{fig:hdamp_Lbar_Alu_HC}
\end{figure}

\Cref{fig:hdamp_Lbar_Alu_HC} presents the dispersion curves obtained with the HC framework. Since no analytical or reference solutions are available for the viscoelastic case, the results are evaluated against expected physical behaviour based on the theoretical insights developed in the preceding examples. Given the absence of geometric symmetry, it is anticipated that mode veering dominates the real wavenumber plots at both the elastic ($s=0$) and viscoelastic ($s=1$) states, with no true crossings occurring.

At the elastic state ($s=0$), uniform wavenumber sampling ($\Delta k = 0.1$) leads to an erroneous crossing (marked by red squares) in the frequency–wavenumber diagram (\cref{fig:hdamp_Lbar_Alu_HC}(a)). This artefact is corrected by the adaptive refinement strategy, which resolves the veering by inserting additional sampling points in regions of rapid eigenvector variation (\cref{fig:hdamp_Lbar_Alu_HC}(b)). The theoretical foundation and numerical validation of this adaptive procedure are detailed in our previous work \cite{xiao_rigorous_2026}.

At the target viscoelastic state ($s=1$), the real wavenumber curves in \cref{fig:hdamp_Lbar_Alu_HC}(c) display clear veering: the two branches approach each other but maintain a positive gap and never cross. This behaviour is consistent with a Type I exceptional point topology. The corresponding imaginary wavenumber curves (\cref{fig:hdamp_Lbar_Alu_HC}(d)) cross over a relatively broad frequency interval, rather than being compressed into an extremely narrow band. The velocity exchange (\cref{fig:hdamp_Lbar_Alu_HC}(f)) likewise unfolds over a finite frequency range, without the sharp, nearly vertical transition that would signal the proximity of an exceptional point to the real axis. These features—a moderate veering gap, a broad imaginary crossing, and a gradual velocity exchange—are hallmarks of a system that remains well within the Type I regime, with the two exceptional points located on opposite sides of, and sufficiently far from, the real axis. The results therefore fully align with the theoretical expectation for a geometry lacking symmetry protection.

A noteworthy feature of this case is the behaviour of the spectral group velocity $v_g$, defined in \cref{eq:group_velocity}, and the energy flux velocity $v_e$, defined in \cref{eq:energy_velocity}. As shown in \cref{fig:hdamp_Lbar_Alu_HC}(f), a slight but discernible offset between the two velocities is observed across the entire frequency range: $v_g$ is consistently slightly larger than $v_e$. This offset is a direct consequence of the non-Hermitian perturbation introduced by the damping, with loss factors $\eta_{\lambda}=10^{-4}$ and $\eta_{\mu}=10^{-3}$. In a purely elastic (Hermitian) waveguide, the left and right eigenvectors satisfy $\mathbf{q}_L = \mathbf{q}^\dagger$, the eigenmodes are orthogonal with respect to the energy inner product, and the two velocities coincide exactly. The introduction of even this modest level of damping breaks the Hermiticity of the system, causing $\mathbf{q}_L \neq \mathbf{q}^\dagger$ and giving rise to the observed offset between the spectral quantity $v_g$ and the physical quantity $v_e$, as discussed in \cref{app:terminology}.

Crucially, despite this offset, the two velocities retain identical structural shapes in the mode-interaction regions: the veering pattern, the locations of the exchange, and the overall curvature are preserved between $v_g$ and $v_e$. This structural congruence confirms that the target viscoelastic system remains well within the Type I regime. At an exceptional point, $v_g$ can exhibit singular behaviour, whereas $v_e$, being derived from integrated physical fields, remains bounded and evolves continuously. The absence of any such singular feature or qualitative shape divergence in the present results indicates that the exceptional points lie far from the real frequency axis, consistent with the broad imaginary crossing and moderate veering gap noted above.

This observation serves a further purpose: it establishes a reference for the behaviour of the $v_g$-$v_e$ discrepancy as a diagnostic tool. In the weakly damped, deeply Type I regime, the two velocities differ by only a small, approximately constant offset and follow the same structural pattern. As damping increases and exceptional points migrate closer to the real axis, the spectral group velocity develops increasingly sharp, localised features that have no counterpart in the energy flux velocity—a signature exploited in \cref{sec:robustness_violated} to identify potential Type II transitions. The present L-bar case therefore provides a baseline against which such anomalous behaviour can be contrasted.

Overall, the HC framework successfully captures the nuanced modal interactions dictated by the underlying exceptional point physics, and the results are fully consistent with the Type I behaviour observed in the previous examples. This example thus confirms the applicability of the proposed method to waveguides of arbitrary cross-section.

\section{Discussion and Analysis}
\label{sec:discussion}
The numerical examples presented in \cref{sec:validation} demonstrate the effectiveness of the proposed homotopy continuation (HC) framework for robust dispersion curve computation across a diverse set of viscoelastic waveguides. In this section, we discuss the method's computational efficiency, analyse its scalability, reliability, and robustness, investigate how it behaves when Type I assumptions are violated (high damping), and compare it with existing approaches in the summary.

\subsection{Computational efficiency}
A primary innovation of the proposed HC framework is the decoupling of the viscoelastic (non‑Hermitian, $s=1$) dispersion problem into two stages: (i) dispersion computation (including mode tracking) for the elastic (Hermitian, $s=0$) waveguide, and (ii) adaptive homotopy continuation from $s=0$ to $s=1$. This strategy significantly reduces the overall computational cost compared to directly solving and tracking modes in the non‑Hermitian system. The efficiency gain arises from three aspects: the use of highly efficient Hermitian eigensolvers in Stage 1, the adaptive sparse homotopy continuation in Stage 2, and parallel processing employed in both stages.

\Cref{tab:efficiency} summarizes the degrees of freedom (DOFs), number of sweeping wavenumbers, number of solution points at $s=0$ and $s=1$, eigengap at mode veering, minimum initial step size, and computational elapsed time for all examples. For all 16‑ply laminates, each layer was discretized into two GLL elements of order 5, yielding 483 DOFs. (Further reduction of the element order, e.g., from 5 to 4, typically preserves accuracy while noticeably reducing the computational cost—a practical option for exploratory parametric studies.)
\begin{table}[ht]
\centering
\caption{Summary of the degrees of freedom (DOFs), number of sweeping wavenumbers, number of solution points at $s=0$ and $s=1$, eigengap at mode veering, minimum initial step size, and computational elapsed time for all examples.}
\label{tab:efficiency}
\footnotesize{
\begin{tabular}{llllllllll}
\toprule
Examples & DOFs & \begin{tabular}{@{}l@{}}Sweeping \\ wavenumbers \\ ($s=0$) \end{tabular} & \begin{tabular}{@{}l@{}}Solution \\ points \\ ($s=0$)\end{tabular} & \begin{tabular}{@{}l@{}}Solution \\ points \\ ($s=1$) \end{tabular} & \begin{tabular}{@{}l@{}}Eigengap \\ at mode \\ veering \end{tabular} & \begin{tabular}{@{}l@{}}Minimum \\ initial \\ step size \end{tabular} & \begin{tabular}{@{}l@{}}Elapsed \\ time \\ (Stage 1) \end{tabular} & \begin{tabular}{@{}l@{}}Elapsed \\ time \\ (Stage 2) \end{tabular} & \begin{tabular}{@{}l@{}}Elapsed \\ time \\ (DC) \end{tabular} \\
\midrule
Sym1    & 483 & 101 & 1146 & 473 & 2.52e-2 & 2.52e-3 & 44.88 s & 86.20 s & 135 s \\
Sym2    & 483 & 100  & 1089 & 443 & 8.67e-2 & 8.67e-3 & 48.13 s & 79.22 s & 223 s \\
UnSym1  & 483 & 272 & 3475 & 968 & 6.22e-3 & 1.00e-3 & 118.47 s & 218.22 s & 280 s \\
UnSym2  & 483 & 177 & 2092 & 806 & 1.01e-2 & 1.01e-3 & 76.21 s & 205.16 s & 287 s \\
Lbar    & 981 & 117 & 1794 & 679 & 4.29e-4 & 1.00e-3 & 285.44 s & 1246.67 s & — \\
\bottomrule
\end{tabular}}
\end{table}

Symmetric laminates (Sym1, and Sym2) have fewer veering regions and a larger eigengap within those regions than unsymmetric laminates (UnSym1 and UnSym2). Consequently, fewer additional sampling points are required during the adaptive wavenumber resampling in Stage 1, leading to a reduced total number of solution points. In fact, the final number of solution points for symmetric laminates is approximately half of that for unsymmetric laminates.

After applying sparse filtering based on the MAC threshold $\bar{\zeta}$ and the interpolation error threshold $\bar{\gamma}$, approximately 40 $\%$ of the solution points are retained for homotopy continuation. Moreover, symmetric laminates generally exhibit larger eigengaps in the veering regions, which allows for a larger initial step size $\Delta s_{\text{init}}$ in the homotopy continuation process. These factors collectively contribute to improved computational efficiency for symmetric laminates.

All computations were performed on a desktop equipped with an Intel i7‑10700 processor (8 cores, 16 threads, base frequency 2.90 GHz) in a Python environment. Parallel processing was employed with 8 concurrent jobs ($n_{\text{jobs}} = 8$) for both the solution of the Hermitian eigenvalue problems in Stage 1 and the homotopy continuation in Stage 2. On average, solving a single Hermitian eigenvalue problem at a given wavenumber requires approximately 0.4–0.5 s for systems with 483 DOFs. This computational cost is influenced by the number of adaptive resampling iterations, which are inherently sequential due to their dependence on previously computed solutions. As symmetric laminates exhibit larger eigengaps, fewer resampling iterations are required, resulting in reduced computational time compared to unsymmetric cases.

A similar trend is observed in Stage 2. The reduced number of solution points, fewer veering regions, and larger initial step sizes in symmetric laminates lead to lower computational costs. For example, for 100 solution points, the elapsed time is 18.22 s (Sym1), and 17.88 s (Sym2), compared to 22.54 s (UnSym1) and 25.45 s (UnSym2). Given the moderate DOFs (483), this level of efficiency is considered satisfactory. 

The L‑shaped bar example, with a larger system size of 981 DOFs, exhibits a similar number of sweeping wavenumbers and solution points as UnSym1. However, due to the significantly smaller eigengap in the veering regions ($4.29\times10^{-4}$), the homotopy continuation becomes substantially more expensive. In this case, homotopy path tracking of 100 solution points requires 178.8 s in Stage 2, indicating that the eigengap is a critical factor governing computational cost.

For comparison, the DC was implemented in MATLAB on the same hardware platform. The results show that, for symmetric laminates Sym1, the total elapsed time of the proposed HC framework (Stages 1+2) is slightly lower than that of DC, while both methods achieve comparable accuracy. For Sym2 (higher damping), HC is considerably faster (127 s vs. 223 s) because DC spends significant extra time in the veering regions, yet still fails to detect all modes and exhibits mode jumps. For the unsymmetric laminates, HC’s total time is similar to or slightly lower than DC’s. However, DC’s runtime for UnSym1 appears lower because it prematurely terminates tracking paths, skipping the most challenging veering regions; this apparent efficiency comes at the cost of reduced accuracy and robustness, as DC fails to capture all modes and produces mode misidentifications. In contrast, the proposed HC framework consistently delivers correct, fully traced dispersion curves, achieving a favorable balance between efficiency and reliability.

Thus, the HC framework provides a robust and efficient alternative for dispersion curve computation in viscoelastic waveguides, especially when high damping or complex modal interactions are present. The practical option of reducing the GLL element order (e.g., from 5 to 4) can further accelerate the overall computation with negligible loss of accuracy, making the method well suited for parametric studies and large‑scale applications.

\subsection{Scalability and applicability}

The scalability of the proposed HC framework to laminates with respect to problem size is investigated using the UnSym1 laminate by systematically varying the thickness discretisation. The laminate is discretised using Gauss–Lobatto–Legendre (GLL) elements with GLL quadrature, a common choice for spectral element methods. A total of 16 configurations are considered, corresponding to combinations of the number of elements per layer ($h=1$–$4$) and element order ($p=3$–$6$), resulting in DOFs ranging from 147 to 1155.

\begin{figure}[htb]
\centering
\includegraphics[width=0.8\columnwidth]{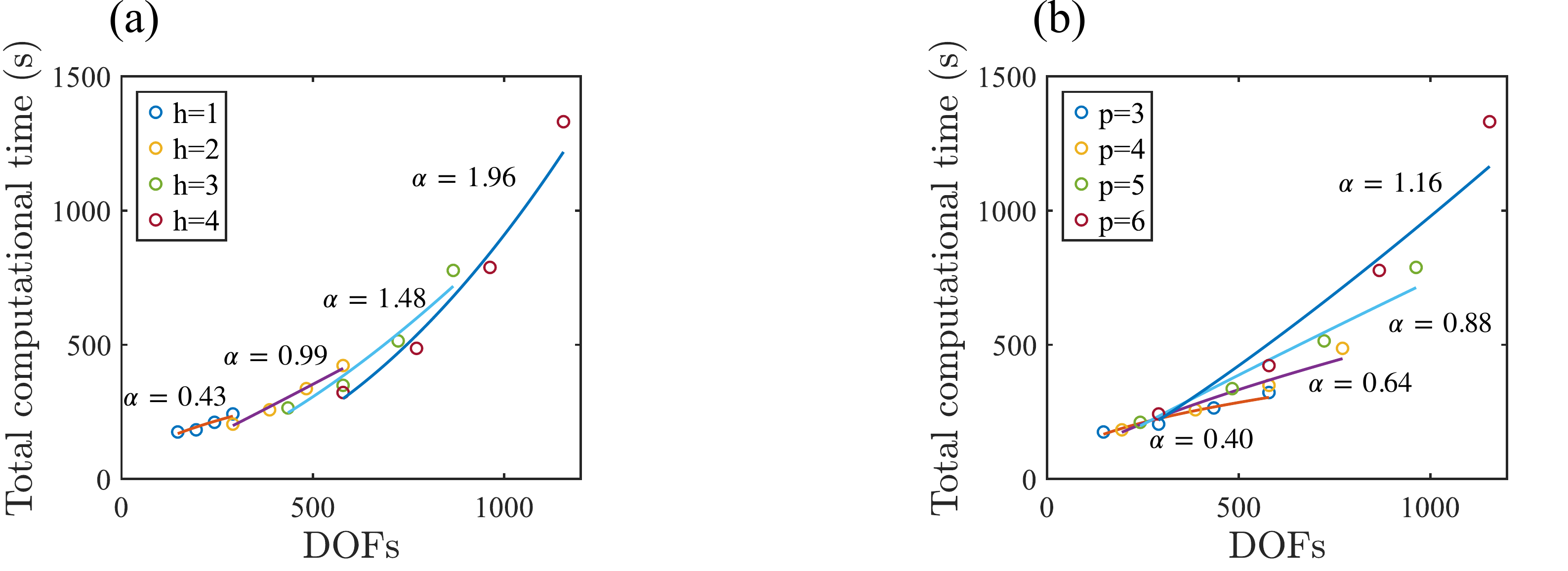}
\caption{Scalability analysis of the HC framework for the UnSym1 laminate. (a) $p$-refinement (fixed number of elements per layer) shows progressively increasing scaling exponents with DOFs. (b) $h$-refinement (fixed element order) exhibits near‑linear scaling, indicating more favourable computational efficiency.}
\label{fig:scalability}
\end{figure}

Computational performance is evaluated in terms of total elapsed time, and the results are presented in \cref{fig:scalability}. To isolate the effects of different discretisation strategies, two complementary analyses are conducted: (i) $p$-refinement at fixed $h$, shown in \cref{fig:scalability}(a), and (ii) $h$-refinement at fixed $p$, shown in \cref{fig:scalability}(b).

For the $p$-refinement cases, the fitted scaling exponents increase with $h$, ranging from $\alpha = 0.43$ to $1.96$. This indicates a transition from sub‑linear to nearly quadratic scaling as the discretisation becomes denser. This behaviour cannot be explained by the increase in DOFs alone; rather, it reflects the growing computational complexity per degree of freedom associated with higher‑order discretisations. As the element order increases, the GLL quadrature leads to denser system matrices and higher condition numbers, while the modal density in the high‑frequency range also increases. These factors collectively render eigenvalue solves more expensive and make homotopy continuation more demanding, often requiring smaller step sizes and a larger number of continuation steps.

In contrast, the $h$-refinement results exhibit significantly more favourable scaling behaviour, with exponents ranging from $\alpha = 0.40$ to $1.16$. In this case, the computational cost grows approximately linearly with DOFs, especially for moderate element orders. This is because increasing the number of elements primarily enlarges the system size while preserving the numerical properties of the local operators (e.g., sparsity and conditioning). Consequently, the cost per degree of freedom remains roughly constant, leading to near‑linear scaling.

The comparison between \cref{fig:scalability}(a) and (b) highlights that the scalability of the proposed HC framework depends not only on the total number of DOFs but also critically on the discretisation strategy. In particular, $p$-refinement increases the intrinsic difficulty of both the eigenvalue solution and the homotopy continuation, whereas $h$-refinement mainly affects problem size without substantially altering numerical conditioning or modal complexity. 

It should be noted that the above analysis is based on laminate thickness discretisation using GLL elements. For two‑dimensional cross‑sections (e.g., the L‑shaped bar), different element types (e.g., quadrilateral or triangular) and integration rules (e.g., standard Gaussian quadrature) are employed. Such choices may lead to different scaling exponents due to variations in matrix sparsity, conditioning, and modal density. Nevertheless, the general trend that $h$-refinement offers better scalability than $p$-refinement is expected to hold across different discretisation strategies.

Overall, the proposed HC framework demonstrates favourable scalability for practical applications. The absence of exponential growth in computational cost confirms that the method remains tractable for large‑scale problems. From an application perspective, $h$-refinement is preferable for improving scalability, while $p$-refinement should be used judiciously when higher accuracy is required, due to its greater impact on computational complexity.

\subsection{Behavior under violated assumptions: entering Type II regime}
\label{sec:robustness_violated}
To examine the numerical behavior of the HC framework when the Type I assumption is violated, we deliberately construct an unsymmetric laminate, UnSym3, with the same stacking sequence as UnSym1 ($[0, 90, 45, -45]_{4}$) but composed of Castaings laminae (\cref{tab:stiffness_matrices}) with the loss modulus uniformly scaled such that the effective loss factor reaches $\eta = 0.05$. This value substantially exceeds the typical range for conventional CFRP materials \cite{gong_improving_2022} and is expected to approach or exceed the critical damping threshold identified in \cref{sec:theory_guarantee}, pushing the system into the Type II regime where EPs migrate to the same side of the real frequency axis.

\begin{figure}[htb]
\centering
\includegraphics[width=1.02\columnwidth]{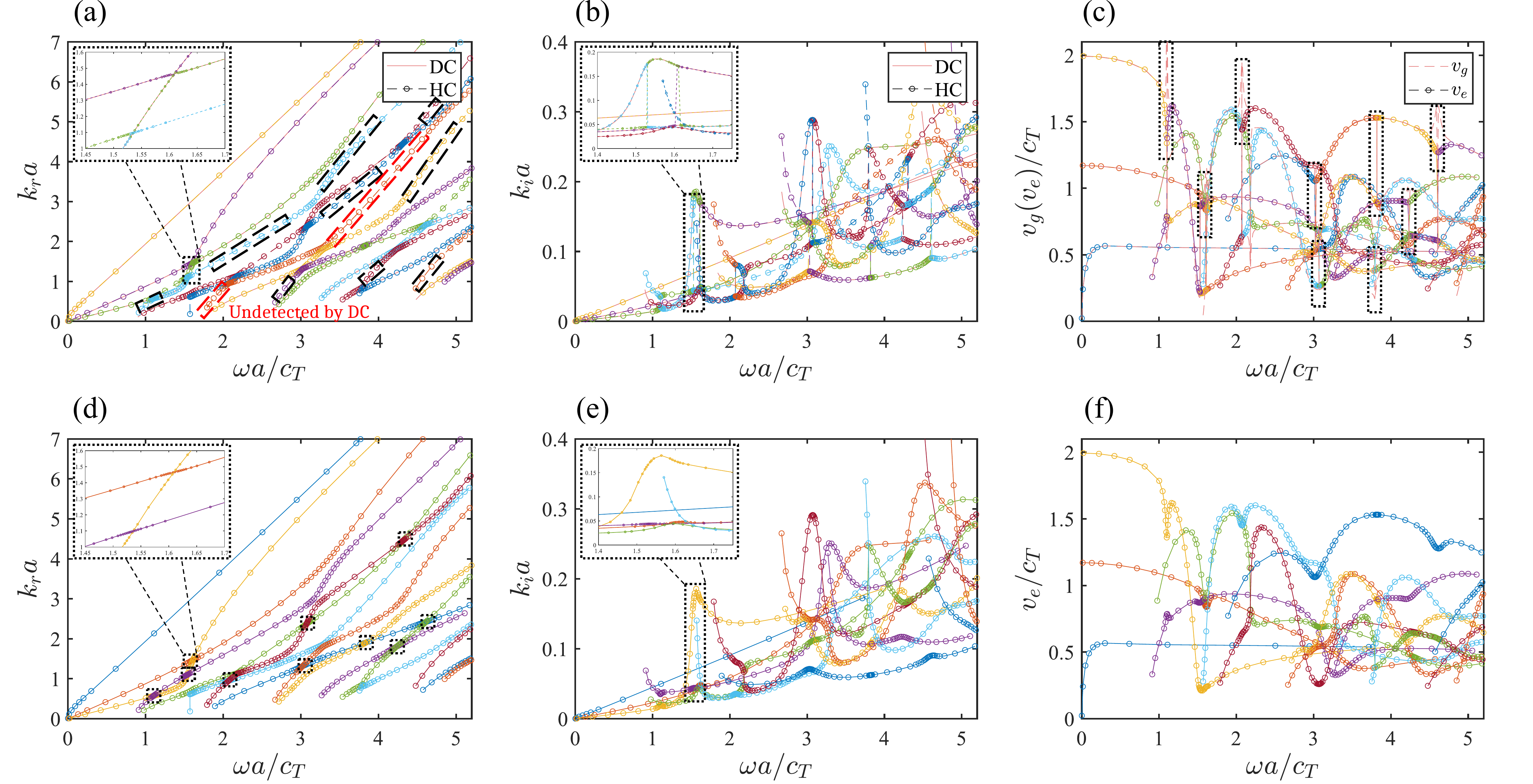}
\caption{Dispersion curves for unsymmetric laminate UnSym3 with an elevated loss factor $\eta = 0.05$. HC results (circles) are superimposed on DC results (solid lines), although the DC method failed to identity all modes. (a) Real wavenumber; (b) imaginary wavenumber; (c) spectral group velocity and energy flux velocity from the HC method, with black squares highlighting veering/crossing regions where a visible discrepancy between the two velocities emerges; (d)–(f) the same curves after post-processing via empirical label exchange: real wavenumber, imaginary wavenumber, and energy flux velocity, respectively. The post-processed curves exhibit real-wavenumber crossings and improved modal continuity, consistent with Type II behaviour.}
\label{fig:violated_assumption}
\end{figure}

\cref{fig:violated_assumption} compares the HC results (circles) with DC (solid lines). The HC solutions coincide with DC wherever the latter successfully tracks modes, confirming that the homotopy continuation faithfully maps elastic solutions to the viscoelastic state even under extreme damping. However, the modal interaction patterns differ markedly from the Type I cases examined in \cref{sec:example_Sym1,sec:example_UnSym1,sec:example_UnSym2,sec:example_Sym2}. The real parts of the wavenumbers exhibit an extremely small veering gap—almost touching—while the imaginary parts cross over a very narrow frequency interval, appearing as nearly vertical lines (\cref{fig:violated_assumption}(a)–(b)). The spectral group velocity $v_g$ and the energy flux velocity $v_e$ diverge visibly in the mode-exchange regions (black squares, \cref{fig:violated_assumption}(c)):
$v_g$ undergoes sharp, localized variations reflecting the rapid change of $dk/d\omega$ along a single Riemann sheet, whereas $v_e$—derived from integrated physical fields—remains bounded and evolves continuously.

These two features—the compressed imaginary crossing and the $v_g-v_e$ discrepancy—match the Type II diagnostic signatures anticipated in \cref{sec:theory_guarantee}. They provide an empirical basis for post-hoc label calibration. Using these indicators, ten veering interactions in the original HC output are identified as likely Type II crossings and relabeled accordingly. \cref{fig:violated_assumption}(d)–(f) present the recalibrated curves: the real wavenumbers now exhibit clear crossings, and the overall modal continuity is visibly improved. It must be emphasized that this recalibration is a posteriori and empirical; the indicators serve as practical warning signals rather than quantitative classifiers of EP topology.

Importantly, the HC continuation itself does not fail in this regime. Because no EP is encountered on the real homotopy path $s \in [0,1]$ in any of the configurations examined—a condition monitored by the Jacobian condition number and adaptive step-size safeguard (\cref{sec:theory_guarantee}), the branch identity along the $s$-axis is preserved, and the numerical wavenumber values are accurate throughout; only the physical mode labels in the observation space require exchange after the transition. The condition number of the Jacobian, monitored during continuation, provides additional diagnostic information if the path approaches an EP in the complex $s$-plane. This demonstrates that the framework provides both numerical fidelity and diagnostic transparency even when the target system lies outside the Type I domain.

\subsection{Robustness and reliability}
The numerical benchmarks in \cref{sec:validation} reveal a systematic difference in failure modes between the HC framework and the Dispersion Calculator (DC), rooted in their fundamentally different algorithmic architectures. DC employs a mode-by-mode frequency-axis continuation strategy: it first detects modes at a high phase velocity cutoff, classifies them into mode families (symmetric, antisymmetric, nonsymmetric Lamb, SH, and Scholte waves), and then traces each mode individually along the frequency axis using predictor–corrector steps at a prescribed frequency resolution. To prevent indefinite stalling when the algorithm loses track of the intended branch—which occurs readily in non-Hermitian veering regions where the Jacobian becomes ill-conditioned—DC employs an internal termination criterion that aborts the current mode and proceeds to the next one.

This architecture creates two distinct failure pathways that are evident in the numerical examples. First, if a mode is not detected during the initial high-phase-velocity sweep—due to a cutoff frequency lying outside the search window or an atypical dispersion slope—it is never traced at all. This explains the entirely missing mode 4 in UnSym1 (\cref{fig:hdamp_Unsym_exa1_CEcast_DCvsHC}(a)) and the two undetected modes in Sym2 (\cref{fig:hdamp_Sym_exa2_DCvsHC}(a), black squares). Second, even when a mode is successfully initialized, the predictor–corrector may lose the branch in a dense veering region and trigger premature termination. In UnSym1, mode 3 jumps to mode 4 near normalized frequency $\omega a/c_T =1.4$ and terminates; mode 7 jumps to mode 3 near 1.9; mode 9 jumps to mode 7 near 2.8 (\cref{fig:hdamp_Unsym_exa1_CEcast_DCvsHC}(a)). In UnSym2, mode 11 terminates at normalized frequency $\omega a/c_T=3.43$ after a significant branch jump, yielding only 30 solution points (\cref{fig:hdamp_Unsym_exa2_CEhern_DCvsHC}(c)). These failures are structurally irreversible: because DC traces modes sequentially and independently, an error at one frequency propagates forward indefinitely; there is no mechanism to recover the lost branch or to verify modal continuity against a reference state.

The HC framework avoids both failure pathways by construction. Mode identification is performed once and for all in the Hermitian elastic regime ($s=0$), where iterative eigensolvers capture the complete spectrum at each wavenumber and the MAC—supported by adaptive refinement—establishes unambiguous connectivity across the entire frequency band (\cref{sec:lossless_stage}). The subsequent homotopy continuation propagates these pre-identified modes along independent paths in the material-parameter space ($s$), which is decoupled from the complex modal interactions in the observation space ($\omega$). Because the continuation parameter is the material loss $s$ rather than frequency $\omega$, veering regions in the $\omega$-plane—however dense or closely spaced—do not cause path divergence in the $s$-plane. The adaptive step-size control (\cref{sec:material_homotopy}) provides an additional safeguard: when local curvature increases as EPs approach the real axis, the algorithm reduces $\Delta s$ automatically rather than jumping branches.

This structural robustness persists even in the Type II regime. As demonstrated in \cref{sec:robustness_violated}, HC continues to produce numerically accurate wavenumbers for UnSym3 ($\eta=0.05$) where DC suffers exacerbated mode loss and requires significantly longer computation time (1476 s versus 616 s for HC). Rather than failing silently or producing unphysical crossings, the framework provides a continuous solution accompanied by diagnostic indicators—an extremely sharp imaginary-part crossing and a discernible $v_g-v_e$ discrepancy—that signal when label exchange is required. The combination of complete spectrum capture at the Hermitian stage, decoupled path tracking in the material space, and diagnostic transparency at the target state distinguishes the HC framework from sequential mode-by-mode approaches that offer no mechanism to detect or correct mode misidentification.

\subsection{Impact of damping-induced EP proximity on computational performance}
\label{sec:eigengap_performance}
While the theoretical guarantee of branch identity continuity holds for all EP-free homotopy paths, the computational performance of the HC framework is governed by the curvature of the homotopy path in the material-parameter space, which in turn reflects the strength of modal interactions in the observation space. Two distinct mechanisms contribute to path curvature, and their distinction is essential for understanding the performance trends in \cref{tab:efficiency}.

First, the loss factor $\eta$ controls the total path length: as $\eta$ increases from 0.003 (Hernando lamina) to 0.02 (Castaings) to 0.05 (UnSym3), the transformation from elastic to viscoelastic becomes more substantial, requiring more continuation steps regardless of local spectral features. This is a global effect: even if the $\omega$-plane EP topology remains Type I with well-separated real parts, a larger $\eta$ demands finer discretization of the $s$-axis.

Second, and more critically, the proximity of $\omega$-plane EPs to the real axis induces local path curvature in the $s$-space through the following mechanism. When an EP approaches the real frequency axis, the two interacting modes exhibit near-coalescence in $k$ at the fixed real $\omega$ used for homotopy continuation. This near-degeneracy renders the dynamic stiffness matrix $\mathbf{D}(k,s)$ ill-conditioned with respect to $k$-perturbations, which in turn degrades the convergence of the Newton corrector in the predictor–corrector algorithm. The adaptive step-size control responds by reducing $\Delta s$, but the path remains traversable because the $s$-plane Jacobian $\partial \mathbf{G}/\partial \mathbf{y} $ (\cref{eq:G_y}) retains full rank as long as no $s$-plane EP is encountered.

The elastic eigengap $\Delta\lambda_{s=0}$ serves as a practical a priori indicator that correlates with both mechanisms. A large elastic eigengap implies that the $\omega$-plane EPs are initially far from the real axis (in the elastic limit, EPs form a symmetric pair whose imaginary separation is proportional to the eigengap). Consequently, even with moderate damping increases, these EPs remain distant, yielding smooth homotopy paths and allowing larger step sizes. This explains the performance contrast between Sym2 ($\Delta\lambda_{s=0} =7.78\times10^{-2}$, 79.22 s Stage 2) and UnSym1 $\Delta\lambda_{s=0} =6.22\times10^{-3}$, 218.22 s Stage 2): the latter's small eigengap anticipates stronger $\omega$-plane modal interactions and thus higher path curvature, despite both laminates operating within the Type I regime.

It is important to emphasize that this performance variation does not compromise robustness. Branch identity continuity depends solely on the absence of $s$-plane EPs on $[0,1]$, not on $\omega$-plane EP proximity. The framework guarantees a one-to-one mapping regardless of path curvature; the curvature only affects the cost of traversal (number of steps, Newton iterations), not the feasibility. For cases approaching the Type II threshold (e.g., UnSym3 with $\eta=0.05$ ), the vanishing $\omega$-plane eigengap signals that physical labels may require post-hoc exchange, but the numerical values remain accurate because the $s$-plane path itself remains regular.

\subsection{Comparison with existing methods}
\Cref{tab:methods_comparison} contrasts the structural paradigms of conventional methods and the proposed HC framework. Conventional approaches—whether direct linearization, iterative eigensolvers, or frequency-axis continuation—attack the non-Hermitian problem at the target state, inextricably coupling eigenvalue solution with mode tracking. The HC framework decouples these tasks by first solving the Hermitian elastic problem and then mapping solutions via homotopy continuation.

\begin{table}[ht]
\centering
\caption{Comparison of existing dispersion calculation approaches with the proposed homotopy continuation framework.}
\label{tab:methods_comparison}
\footnotesize{
\begin{tabular}{lll}
\toprule
Feature / Aspect
& Conventional Approaches
& Proposed HC Framework
\\
\midrule
Problem-Solving Paradigm
& \begin{tabular}{@{}l@{}} Directly solves non-Hermitian problem \\ at target state \end{tabular} 
& \begin{tabular}{@{}l@{}} Decouples into Hermitian stage + \\ homotopy tracking \end{tabular} 
\\
Mode Tracking
& \begin{tabular}{@{}l@{}}  Post-processing heuristics (MAC, sorting) \\ at target state \end{tabular} 
& \begin{tabular}{@{}l@{}} Established at lossless stage, \\ theoretically preserved \end{tabular} 
\\
Theoretical Guarantee
& None 
& \begin{tabular}{@{}l@{}}Rigorous guarantee of branch identity \\ continuity along homotopy path \end{tabular} \\
Handling of Mode Veering
& \begin{tabular}{@{}l@{}} Fragile; may miss modes or produce\\ artificial jumps \end{tabular} 
& \begin{tabular}{@{}l@{}}Robust; path tracking follows \\ correct branch \end{tabular} \\
\begin{tabular}{@{}l@{}} Extensibility to \\ Frequency-Dependent Damping \end{tabular} 
& \begin{tabular}{@{}l@{}} Possible but challenging; increased solution \\ complexity and tracking difficulty \end{tabular}
& \begin{tabular}{@{}l@{}} Straightforward under Type I topology \\ (analytic parameter dependence) \end{tabular}  \\
Computational Efficiency 
& Parameter-sensitive; conservative settings
& \begin{tabular}{@{}l@{}} Efficient due to sparse mapping \\ and parallelizability \end{tabular}  \\
\bottomrule
\end{tabular}}
\end{table}

This structural difference manifests in both robustness and efficiency. Because HC performs mode tracking in the Hermitian domain, it captures the full spectrum without heuristic shift selection or manual initialization. The sparse key-point mapping (retaining approximately 40$\%$ of dense solutions, \cref{sec:lossless_stage}) and parallelizable independent paths further reduce computational cost. \Cref{tab:efficiency} shows that for symmetric laminates with moderate damping (Sym1, Sym2), HC achieves total times comparable to or lower than DC, while reliably detecting all modes. For unsymmetric laminates with pervasive veering (UnSym1, UnSym2), HC remains efficient despite smaller eigengaps, whereas DC prematurely terminates branches to avoid veering regions—an apparent time saving that actually reflects incomplete output.

The progressive nature of the homotopy transformation confers an additional advantage for the weakly damped materials that dominate practical engineering composites. Because each predictor–corrector step advances the solution by a small material increment $\Delta s$ , the total path length is short when the target loss factor is small (e.g., $\eta \approx 0.003$  for Hernando lamina). In this regime, the transformation from elastic to viscoelastic is numerically trivial, and the framework operates effectively as a perturbed elastic solver with guaranteed modal connectivity.

At elevated damping ($\eta = 0.05$), the efficiency gap widens further. DC requires significantly longer computation time as damping increases, because dense modal interactions force conservative step sizes and repeated re-initialization. HC, by contrast, handles the increased path curvature through adaptive step-size control without sacrificing solution completeness. The numerical comparison substantiates this: for UnSym3, DC requires 1476 s compared to HC's 616 s, yet still produces discontinuous curves with missing branches.

The framework is theoretically extensible to frequency-dependent damping models (Kelvin–Voigt, Zener) when the Type I topology holds, because these models satisfy the same analyticity condition that underpins the branch identity guarantee. The current implementation focuses on hysteretic damping due to the availability of reference tools (DC) for validation; extension to frequency-dependent models is planned for release. Overall, the HC framework offers a favorable trade-off: it provides the theoretical guarantees of branch identity (\cref{sec:theory_guarantee}) while delivering superior numerical robustness and competitive—or better—computational efficiency across the damping range tested.

\section{Conclusion and future work}
\label{sec:conclusion}
This paper has demonstrated that the fundamental challenge of non-Hermitian dispersion analysis—the inseparability of eigenvalue solution and mode tracking—can be resolved by an inter-manifold transport strategy. Rather than confronting the lossy non-Hermitian problem directly, the framework establishes globally unique modal identities on the Hermitian elastic anchor and transports them to the viscoelastic target along a material homotopy path. Analytic perturbation theory guarantees that branch identity is preserved throughout this transport for all EP-free paths, irrespective of the target state's topology. Yet numerical continuity does not suffice for physical correctness: whether the transported labels are inherited automatically or require post-hoc exchange is governed by the Type I/II exceptional-point topology in the observation space. Type I configurations (EPs on opposite sides of the real frequency axis) permit seamless label inheritance, whereas Type II configurations (an EP crossed the real axis) demand a swap to restore physical consistency. This separation of numerical continuity (guaranteed by the homotopy) from physical label validity (determined by EP topology) clarifies why existing methods fail—by attacking the non-Hermitian problem directly, they irreversibly couple numerical error with topological ambiguity—and provides a transparent, theoretically grounded alternative.

Numerical benchmarks across symmetric laminates, unsymmetric laminates with pervasive veering, and an L-shaped bar confirm that the framework produces physically consistent dispersion curves in the Type I regime (validated for $\eta \approx 0.003$), with automatic label inheritance. For the Type II regime ($\eta=0.05$), the framework continues to yield numerically accurate wavenumbers but requires post-hoc label exchange, guided by two diagnostic signatures: extremely sharp imaginary-part crossings and spectral–energetic velocity discrepancies. These diagnostics transform the framework from a ``black-box'' solver into a transparent tool that alerts users when assumptions are violated.

The main conclusions of this work are as follows:
\begin{itemize}
\item \textbf{Inter-manifold transport guarantees branch identity continuity irrespective of target-state topology.} The framework successfully decouples mode tracking from non-Hermitian eigensolving by establishing modal identities on the Hermitian elastic anchor and transporting them along a material homotopy path. Analytic perturbation theory ensures a one-to-one correspondence between $s=0$ and $s=1$  for all EP-free paths, while the separation of numerical continuity from physical label validity means that even when Type II configurations violate automatic inheritance, the underlying dispersion solutions remain numerically valid and require only post-hoc label adjustment.

\item \textbf{Non-Hermitian exceptional-point topology governs observable dispersion patterns and label inheritance.} Numerical validation confirms that the Type I/II EP classification—previously unexplored in guided-wave dispersion analysis—directly dictates whether dispersion branches exhibit real-part veering with imaginary-part crossing (Type I) or real-part crossing with imaginary-part veering (Type II). This correspondence provides a physically grounded criterion for label inheritance after homotopy transport: Type I configurations retain anchor labels automatically, whereas Type II configurations require a swap. 

\item \textbf{Validated practical applicability across damping regimes and geometries.} For weakly damped CFRP laminates ($\eta \approx 0.003$), numerical evidence confirms Type I topology across diverse configurations—symmetric, unsymmetric, and arbitrary cross-section—yielding automatic label inheritance without post-processing. At elevated damping ($\eta = 0.05$), the framework maintains numerical correctness where conventional methods fail entirely, resolving dense modal interactions through adaptive step-size control. Two empirical diagnostic signatures—an extremely sharp imaginary-part crossing and a marked spectral--energetic velocity discrepancy—reliably flag Type II transitions when a priori prediction is unavailable.

\item \textbf{Diagnostic transparency for operational deployment.} The framework transforms non-Hermitian dispersion solvers from black-box tools into transparent systems that alert users to assumption violations. When the Type I assumption is violated, the solver returns numerically accurate wavenumbers accompanied by two physically motivated diagnostic signatures that guide label exchange. This transparency is critical for engineering practice, where users require confidence in both the numerical validity of solutions and the physical correctness of assigned mode labels.
\end{itemize}

Despite its demonstrated robustness, the proposed framework has several limitations that warrant further investigation. First, the homotopy method assumes that no EP lies on the real $s$-axis. As the damping level of the target viscoelastic system increases, the range of material states spanned by the homotopy path $s \in [0,1]$ broadens, increasing the likelihood that an EP in the complex $s$-plane approaches the real axis. When this occurs, the Jacobian matrix becomes increasingly ill-conditioned, and the adaptive step-size control responds by progressively reducing. If the path passes too close to an EP, the algorithm terminates prematurely without returning a solution at $s=1$. The condition number of the Jacobian, which is monitored during continuation, can provide diagnostic information in such situations. Consequently, the present framework is not designed for strongly damped regimes where the loss modulus $\mathbf{C}''$ approaches or exceeds the same order of magnitude as the elastic modulus $\mathbf{C}'$, as the extended homotopy path significantly raises the risk of encountering an EP.

Second, while the theoretical guarantee of branch identity continuity holds for any damping model with analytic parameter dependence, the current implementation and numerical validation focus on hysteretic damping. The framework is, in principle, applicable to frequency-dependent models such as Kelvin–Voigt and Zener without modification, because at each fixed frequency the loss factor of such models is a constant value and only the effective length of the homotopy path differs from the hysteretic case. However, the absence of a reference method capable of handling frequency-dependent damping (DC is limited to hysteretic damping) has precluded a systematic validation for these models. The implementation will be made publicly available at \url{https://github.com/dongxiao96/TopoDisper} upon publication.

Third, the framework lacks a priori quantitative criteria to determine whether a given system lies in the Type I or Type II regime; only post‑hoc diagnostics—specifically, the two empirical indicators based on imaginary‑part crossing sharpness and the $v_g$-$v_e$—can currently indicate a potential assumption violation. Importantly, even when the system enters the Type II regime, the homotopy framework still returns numerically accurate dispersion solutions; only the automatic inheritance of physical labels is lost, requiring the post‑hoc label exchange described in \cref{sec:robustness_violated}.

Building on these findings, future work will proceed along the following directions:
\begin{itemize}
\item \textbf{Quantitative prediction of Type II occurrence}: The post‑hoc diagnostics developed in this work provide reliable empirical indicators of EP proximity and potential Type II transitions. Building on these insights, future work will develop explicit EP tracking algorithms in the complex frequency plane to establish a priori quantitative criteria for the Type I/II classification. By locating the critical damping threshold at which EPs cross the real frequency axis, such a tool would complement the present framework with a pre‑computation diagnostic, replacing the current empirical indicators with a rigorous, predictive criterion.

\item \textbf{Advanced continuation for strongly damped regimes}: Explore high‑order continuation methods based on Puiseux series and regularized eigenvalue pairs to overcome the numerical limitations near EPs \cite{nennig_high_2020}. These techniques can reconstruct regularized functions that remain analytic across the singularity, enabling stable path tracking even when EPs lie extremely close to the real $s$-axis. This would extend the applicability of the homotopy framework to the strongly damped materials where the current predictor–corrector algorithm may terminate prematurely, and would also provide an alternative route to computing dispersion curves across Type II transitions that automatically handles the associated label exchange.

\item \textbf{Extension to other non‑Hermitian systems}: The proposed homotopy framework is not limited to guided wave problems; it can be extended to other physical systems governed by non‑Hermitian eigenvalue problems with an analytic parameter transition. Examples include fluid‑loaded structures, acoustic–elastic coupled systems, and structures subject to external damping or gyroscopic effects. Such extensions would demonstrate the broader applicability of the method and its potential to serve as a general‑purpose tool for robust eigen‑tracking across a wide class of non‑Hermitian problems.
\end{itemize}

In summary, the proposed adaptive homotopy continuation framework offers a theoretically grounded, robust, and efficient approach for dispersion curve computation in viscoelastic waveguides. By inter-manifold transport—confining mode tracking to the Hermitian anchor and propagating identities along a material homotopy—it overcomes the fragility that plagues existing methods in regions of veering and strong modal interaction. A rigorous theoretical foundation, grounded in analytic perturbation theory, guarantees branch identity continuity along the homotopy path, while the automatic inheritance of physical mode labels is established for systems retaining a Type I EP topology. For systems that transition to Type II, the framework continues to return numerically accurate solutions; two physically motivated post‑hoc diagnostics alert the user when label exchange may be required. The numerical validation across diverse waveguide configurations and damping levels confirms the method's accuracy, robustness, and practical utility. With the planned extensions, the framework is poised to become a versatile tool for non-Hermitian guided wave analysis in complex engineering structures.

\appendix
\section{Appendix: Challenges in non-hermitian dispersion analysis}
\label[appendix]{app:challenges}
Extracting accurate dispersion curves from viscoelastic waveguides requires solving the non-Hermitian eigenvalue problem presented in \cref{eq:SAFE_quad_eigen} and subsequently tracking the identified modes across frequency to assemble continuous modal branches. This subsection reviews the principal approaches reported in the literature, highlighting their inherent limitations, particularly for materials characterized by hysteretic damping. The core difficulty stems from the non-Hermitian nature of the system itself: existing methods capable of calculating dispersion curves for attenuated modes can be troublesome to implement and the solutions are not as reliable as in the perfectly elastic case \cite{quintanilla_guided_2015}.

\subsection{Discrete frequency methods: linearization and root-searching}
A common strategy linearizes the quadratic eigenvalue problem into a larger generalized eigenvalue problem. By introducing an auxiliary vector, \cref{eq:SAFE_quad_eigen} can be rewritten as:
\begin{equation}
    \begin{aligned}
    & \mathbf{A}(\omega) \mathbf{z}  = k\mathbf{B}\mathbf{z},  \mathbf{z} = [\mathbf{q}^T, \; (k\mathbf{q})^T]^T,  \\
        & \mathbf{A} = \left[ \begin{matrix}
\mathbf{0} & \mathbf{I} \\
 -\mathbf{K}_1 + \omega^2\mathbf{M} & - i\mathbf{K}_{2}
\end{matrix} \right], \; 
       \mathbf{B}  = \left[ \begin{matrix}
\mathbf{I} &  \mathbf{0}  \\
\mathbf{0} &  \mathbf{K}_{3}
\end{matrix} \right].
    \end{aligned}
\end{equation}
which is a linear eigenvalue problem of dimension $2n$  (where $n$ is the original number of degrees of freedom). For small to medium-scale models, a dense direct eigensolver (e.g., \texttt{eig} in MATLAB) can be applied to compute \emph{all} $2n$ eigenvalues at a given frequency. Although straightforward and capable of producing the complete spectrum, this approach scales as $\mathcal{O}(n^3)$ in computation and $\mathcal{O}(n^2)$ in memory, making it impractical for large-scale finite element models.

For larger models, iterative eigensolvers such as the shift-invert Arnoldi method \cite{lehoucq_deflation_1996} are preferred (e.g., \texttt{eigs} in MATLAB). These solvers compute eigenvalues near a user-defined target shift, drastically reducing computational cost by leveraging sparsity. While this approach benefits from mature numerical libraries, it presents several challenges for hysteretic damping. The non-Hermitian nature of the system makes the eigenvalue problem intrinsically more difficult; iterative solvers may exhibit poor convergence, and the choice of shift parameters is often heuristic, risking the omission of physically relevant modes. Unlike Hermitian problems, where a single shift can often capture a complete spectrum within a frequency band, non-Hermitian systems require exploring multiple shift parameters across the complex plane to ensure that no valid solution is missed.

To address the limitations of shift-based iterative solvers, contour integral methods have been proposed as a more robust alternative for extracting all eigenvalues within a specified region of the complex plane \cite{sakurai_projection_2003}. These methods, successfully applied to viscoelastic waveguide problems \cite{mazzotti_coupled_2013}, transform the nonlinear eigenvalue problem into a linear one inside a chosen contour, offering improved reliability without requiring initial guesses. However, they come at increased computational cost and still require exploring different regions of the complex plane to obtain a complete set of modes.

An alternative family of methods circumvents linearization by treating the characteristic equation directly as a transcendental function in the wavenumber $k$. Approaches such as the GMM, SMM, and complex root-searching algorithms \cite{orta_comparative_2022, quiroga_evaluation_2025} search the complex $k$-plane for zeros of the determinant or related functions. While accurate for isolated modes, they become cumbersome in densely populated mode regions or when modes exhibit veering and crossing. The root-searching process must be repeated for each mode at each frequency, and ensuring that all roots are found—without omission or duplication—requires careful initialization and often manual intervention. Like direct eigensolvers, these methods also operate on a frequency-by-frequency basis and offer no inherent modal connectivity.

For both linearization-based and root-searching methods, mode tracking is performed as a post-processing step. In elastic (Hermitian) waveguides, the Modal Assurance Criterion (MAC) \cite{allemang_modal_2003} has proven reliable:
\begin{equation} \label{eq:MAC}
    \begin{aligned}
    \mathrm{MAC}[\mathbf{q}_i, \mathbf{q}_j] = \frac{\left|\mathbf{q}_i^{\dag} \mathbf{M} \mathbf{q}_j\right|^2}{ \left(\mathbf{q}_i^{\dag} \mathbf{M} \mathbf{q}_i\right) \left(\mathbf{q}_j^{\dag} \mathbf{M} \mathbf{q}_j\right)},
    \end{aligned}
\end{equation}
where $\dag$ denotes the Hermitian conjugate, which is widely used in non-Hermitian physics. However, in viscoelastic (non-Hermitian) systems, the loss of orthogonality and the presence of near-coalescing modes degrade MAC reliability. Moreover, these methods operate on a frequency-by-frequency basis, preventing the reuse of information across frequencies and resulting in a discrete set of eigenvalues that must subsequently be associated into continuous modes \cite{mazzotti_25d_2013}.

\subsection{Frequency continuation and mode tracking difficulties}

Numerical continuation methods, particularly arc-length continuation along the frequency axis \cite{allgower_numerical_1990}, offer an efficient alternative for tracing individual modes. Starting from a known eigensolution at an initial frequency $\omega_0$, the dispersion curve for a given mode is traced by solving an extended system that enforces an arc-length constraint \cite{maruyama_continuation_2025}:
\begin{equation} 
    \begin{aligned}
        & \mathbf{F}(k, \omega, \mathbf{q}) = 0, \\
        & \dot{\mathbf{q}}^{\dagger}(\mathbf{q} - \mathbf{q}_{pred}) 
        + \dot{k}^{*}(k-k_{pred})
        + \dot{\omega}(\omega-\omega_{pred}) = 0,
    \end{aligned} 
\end{equation}
where $\mathbf{F}(k, \mathbf{q}, \omega) = [\mathbf{K}_1 + i k \mathbf{K}_2 + k^2 \mathbf{K}_3 - \omega^2 \mathbf{M}] \mathbf{q}$. The predictor-corrector algorithm proceeds with tangential prediction and Newton-based correction, naturally handling turning points in the solution curve. The primary advantage of this approach lies in its efficiency: the number of continuation steps is typically far smaller than the number of discrete frequency points required for a full sweep, and the resulting curves are inherently smooth and continuous.

However, applying this approach directly to the viscoelastic problem (i.e., performing continuation in frequency at the target damping state) introduces fundamental difficulties. First, it requires a high-quality starting solution at some initial frequency—a complex eigensolution that is itself challenging to obtain reliably. Obtaining such starting solutions typically involves solving the non-Hermitian eigenvalue problem multiple times across the complex plane to ensure that all modes have been identified at the starting frequency, effectively reintroducing the difficulties of the direct methods.

Second, and more critically, continuation along the frequency axis is highly sensitive to mode veering and degeneracy. In regions where two modes approach closely, the Jacobian matrix of the extended system becomes ill-conditioned, leading to convergence failure or spurious mode jumping. The non-Hermitian nature of the viscoelastic problem further exacerbates this sensitivity. Unlike Hermitian systems, where mode veering can be reliably tracked with sufficient frequency resolution, non-Hermitian systems may exhibit exceptional points where eigenvalues coalesce and eigenvectors become parallel. At such points, the Jacobian becomes singular and the continuation method fails. Moreover, mode crossings in non-Hermitian systems are not necessarily symmetry-protected and can occur, making it impossible to guarantee that a continuation path will remain on the intended modal branch. These challenges render direct frequency-axis continuation at the target damping state inherently fragile.

Unlike discrete frequency methods, frequency continuation integrates mode tracking into the solution process; it does not rely on post-processing criteria such as MAC. However, its fragility near veering and exceptional point regions remains a critical limitation.

\subsection{Universal mode tracking challenges}
Independent of the numerical strategy employed to obtain eigenvalues, the association of solutions at adjacent frequencies into continuous modal branches—mode tracking—presents its own set of challenges. In elastic waveguides modeled by the SAFE method, the system matrices are Hermitian, and modal behavior is well-understood. Although eigenvectors are complex-valued due to the presence of the $ik\mathbf{K}_2$ term and numerical error, the Hermitian structure imposes strong constraints: mode crossings are symmetry-protected, and mode veering occurs in a predictable manner. The Modal Assurance Criterion (MAC) has proven to be a reliable tool for linking modes across frequencies in such settings.

In viscoelastic waveguides, the loss of Hermiticity removes these structural constraints, making mode tracking substantially more difficult. Mode veering becomes more complicated to analyze and track, as the rapid exchange of eigenvector characteristics can occur over narrower frequency intervals. Mode degeneracy in non-Hermitian systems can manifest as exceptional points, where eigenvalues coalesce and eigenvectors become parallel—a phenomenon without analogue in Hermitian systems. These exceptional points introduce topological complexity that cannot be resolved by simple frequency refinement alone.

Existing mode tracking strategies for non-Hermitian systems typically fall into two categories, each with significant limitations. The first relies on similarity measures such as the MAC. However, the reliability of MAC degrades considerably in non-Hermitian systems due to the loss of orthogonality and the presence of near-coalescing modes, where eigenvectors become nearly parallel and the inner product loses its discriminative power. The second employs physical heuristics—for example, sorting modes by the real part of the wavenumber or by group velocity—but these rules fail in regions of mode veering or crossing, where modal identities can exchange without clear signatures in the sorted order. As a result, mode tracking in the viscoelastic regime often demands substantial manual supervision, especially when damping is significant or the frequency range includes complex modal interactions.

Crucially, mode tracking difficulties are compounded by the underlying eigenvalue solution challenges. If the eigenvalue solver fails to capture all physical modes at a given frequency—a common risk in non-Hermitian problems with poorly chosen shift parameters—then no subsequent tracking algorithm can recover the missing branches. Errors in eigenvalue computation propagate irreversibly into mode tracking, and once modes are misassigned, they cannot be corrected by post-processing alone.

\subsection{Limitations of existing non-Hermitian dispersion methods}
The preceding review exposes a common structural limitation that cuts across all existing methods: they operate directly on the non‑Hermitian target manifold and attempt to track modes using only local information—eigenvector similarity, tangent prediction, or heuristic sorting. In other words, they are confined to \emph{intra‑manifold tracking}. None of these approaches possess any awareness of the global Riemann‑sheet topology on which the eigenvalues reside. When exceptional points (EPs) lie close to the real frequency axis, the eigenvectors of interacting modes become nearly indistinguishable, and the local similarity measures on which trackers rely degrade, causing mode tracking to fail—often silently, by producing numerically smooth curves with physically incorrect labels. This fragility is not a shortcoming of any particular solver or tracking algorithm, but a direct consequence of the topological blindness inherent in intra‑manifold operation. Neither dense direct eigensolvers paired with MAC, nor iterative shift‑invert or contour‑integral methods combined with post‑processing, nor even frequency‑axis continuation can guarantee correct modal identity preservation across the lossless‑to‑lossy transition, because they all lack a global reference for mode identities. The present work replaces this blind intra‑manifold tracking with inter‑manifold transport, where identities are defined once on a well‑posed Hermitian anchor and safely transported to the viscoelastic target via homotopy continuation.

\section{Appendix: Terminology and key concepts in non-Hermitian dispersion analysis}
\label[appendix]{app:terminology}
This section provides precise definitions of key concepts used in the analysis of non-Hermitian dispersion systems governed by the SAFE formulation $\mathbf{D}(k,\omega)\mathbf{q} = \mathbf{0}$ (\cref{eq:SAFE_quad_eigen}), where $k \in \mathbb{C}$ is the wavenumber, $\omega \in \mathbb{C}$ is the angular frequency, and $\mathbf{q}$ is the right eigenvector (for simplicity, $\mathbf{q}$ is used to represent the right eigenvector rather than $\mathbf{q}_R$). In non-Hermitian systems, the dispersion relation is multi-valued and must be interpreted through its analytic structure on a Riemann surface. A clear distinction between analytic objects and physically tracked quantities is therefore essential.

\noindent\textbf{Dispersion relation.}
The dispersion relation is defined implicitly by
\begin{equation}
\det \mathbf{D}(k,\omega) = 0,
\end{equation}
which establishes a multi-valued functional relationship between $k$ and $\omega$. This relation defines a complex algebraic curve whose solutions are naturally interpreted on a Riemann surface.

\noindent\textbf{Riemann surface and Riemann sheets.}
The multi-valued function $k(\omega)$ becomes single-valued when lifted onto a Riemann surface composed of multiple sheets. Each \emph{Riemann sheet} corresponds to one analytic continuation of the dispersion solution. Branch points (such as exceptional points) connect different sheets, while branch cuts define discontinuities that arise when projecting the surface onto the complex plane.

\noindent\textbf{Branches and branch identity.}
A \emph{branch} is a locally single-valued analytic function $k^{(i)}(\omega)$ defined on a Riemann sheet. The \emph{branch identity} labels a specific analytic continuation of the dispersion relation. Under continuous parameter variation that does not pass through a branch point, the branch identity is preserved. However, when restricted to real-frequency slices, apparent discontinuities may arise due to intersections with branch cuts.

\noindent\textbf{Exceptional points}. An exceptional point (EP) is a branch point in the complex parameter plane where two or more eigenvalues and their corresponding eigenvectors simultaneously coalesce, and the matrix becomes defective \cite{kato_perturbation_1995, heiss_exceptional_2004}. For the SAFE system $\mathbf{D}(k,\omega)\mathbf{q} = \mathbf{0}$, an EP ($k_{\text{EP}}, \omega_{\text{EP}}$) is defined by the simultaneous satisfaction of the dispersion relation and the stationarity condition \cite{ghienne_beyond_2020}:
\begin{equation}
\det \mathbf{D}(k,\omega) = 0, \qquad
\frac{\partial}{\partial \omega} \det \mathbf{D}(k,\omega) = 0.
\end{equation}
At an EP, the algebraic multiplicity of the eigenvalue exceeds its geometric multiplicity; the eigenbasis collapses, meaning that two or more originally independent eigenvectors degenerate into a single direction in the state space. This collapse is the mathematical origin of the extreme sensitivity and modal non-orthogonality observed near EPs. For a pair of interacting modes in Hermitian system, EPs appear as a conjugate pair in the complex frequency plane, and their distance from the real axis controls the strength of the associated veering or crossing interaction \cite{heiss_phases_1999}.

\noindent\textbf{Type I and Type II EP topology}. The configuration of an EP pair relative to the real frequency axis determines the observable crossing behaviour along the real frequency axis \cite{keck_unfolding_2003}. In a \emph{Type I} configuration, the two EPs lie on opposite sides of the real axis. Scanning the real frequency axis then yields real-part veering (an avoided crossing with a positive gap) accompanied by imaginary-part crossing. The order of the real parts does not change across the interaction region, and the physical mode shapes associated with each branch remain consistent. In a \emph{Type II} configuration, both EPs lie on the same side of the real axis. The observable behaviour reverses: the real parts cross while the imaginary parts veer. To maintain physical continuity in this case, the mode labels must be exchanged at the crossing point. In the elastic limit ($s=0$), EPs form a complex-conjugate pair symmetric about the real axis, corresponding to a Type I topology. As material damping increases, the EPs migrate continuously in the complex plane; whether a transition to Type II occurs depends on the damping magnitude relative to a problem-dependent critical threshold. This classification is fundamental to the mode identity guarantee developed in \cref{sec:theory_guarantee}.

\noindent\textbf{Dispersion curves.}
Dispersion curves are the restriction of the dispersion relation to real-valued frequency (or wavenumber), typically expressed as $k(\omega)$ with $\omega \in \mathbb{R}$. These curves are projections of the underlying Riemann surface and do not uniquely encode branch identity.

\noindent\textbf{Modes.}
A \emph{mode} is defined as a solution pair $(k(\omega), \mathbf{q}(\omega))$ satisfying $\mathbf{D}(k,\omega)\mathbf{q} = \mathbf{0}$, 
together with the corresponding left eigenvector $\mathbf{q}_L$ satisfying $\mathbf{q}_L^\dagger \mathbf{D}(k,\omega)=\mathbf{0}$. Mathematically, a mode corresponds to an \emph{analytic continuation of an eigenpair on the Riemann surface}. Away from degeneracies, each mode is associated with a single branch. However, near exceptional points, eigenvalues and eigenvectors coalesce, and this one-to-one correspondence breaks down.

\noindent\textbf{Continuity of modes.}
A central concept is that a mode is not merely a pointwise solution, but a \emph{continuous (analytic) trajectory} on the Riemann surface:
\begin{equation}
\omega \;\mapsto\; \left(k(\omega), \mathbf{q}(\omega)\right).
\end{equation}
This continuity is defined in the sense of analytic continuation with respect to $\omega$ (or other parameters). It ensures that the mode represents a physically meaningful evolution of a wave solution, rather than a collection of unrelated eigenpairs.

\noindent\textbf{Mode tracking and physical mode label.}
In practical computations, dispersion solutions are obtained at discrete frequency points, and such continuous trajectories are not directly available. \emph{Mode tracking} is therefore the procedure of reconstructing a sequence of eigenpairs
\begin{equation}
\left\{(k(\omega_i), \mathbf{q}(\omega_i))\right\}_{i=1}^N
\end{equation}
that approximates a continuous path on the Riemann surface by enforcing continuity criteria across successive frequencies.

A \emph{physical mode label} is then assigned to each such tracked continuous trajectory, serving as an identifier of the same underlying mode. Importantly, the label does not define the mode itself; rather, it is attached \emph{after} a continuous path has been established. Physical quantities such as phase velocity, attenuation, symmetry, or energy content may assist in the tracking process or in interpreting the resulting modes, but they do not constitute the defining criterion. The essential requirement is the continuity of the eigenpair sequence, which ensures that the labeled mode represents a physically meaningful and consistently evolving solution.

\noindent\textbf{Spectral group velocity.}
The spectral (dispersion-derived) group velocity is defined as
\begin{equation}
v_g = \left(\frac{dk}{d\omega}\right)^{-1},
\end{equation}
which can be evaluated from the implicit eigenvalue problem as
\begin{equation} \label{eq:group_velocity}
v_g = -\frac{\mathbf{q}_L^\dagger  \partial_k \mathbf{D}  \mathbf{q}}{\mathbf{q}_L^\dagger  \partial_\omega \mathbf{D}  \mathbf{q}}.
\end{equation}
This quantity reflects the local geometry of the dispersion relation on the Riemann surface.

\noindent\textbf{Energy flux velocity.}
The energy flux velocity is defined from physical energy transport as
\begin{equation} \label{eq:energy_velocity}
v_e = \frac{\displaystyle \int_{\Omega} S(y,z) dA}{\displaystyle \int_{\Omega} W(y,z) dA},
\end{equation}
where $S(y,z)$ is the time-averaged energy flux and $W(y,z)$ is the stored energy density. These quantities are obtained from reconstructed physical fields and therefore incorporate spatial structure and modal interactions.

\noindent\textbf{Spectral--energetic distinction.}
In Hermitian systems, the spectral group velocity $v_g$ and the energy flux velocity $v_e$ coincide, owing to the orthogonality of eigenmodes and the consistency of the associated inner products. In non-Hermitian systems, however, the left eigenvectors differ from the Hermitian conjugates of the right eigenvectors ($\mathbf{q}_L \neq \mathbf{q}^\dagger$), introducing modal non-orthogonality and a fundamental discrepancy between the two velocity measures. This discrepancy is strongly amplified near exceptional points (EPs): while the continuity of modes on the Riemann surface is preserved, the equivalence between dispersion-derived and energy-based velocities breaks down. As eigenvectors coalesce, $v_g$ can exhibit increasingly sharp, even singular, mode-exchange behaviour, whereas $v_e$, derived from integrated physical fields, necessarily remains bounded and evolves continuously. Consequently, a growing divergence between $v_g$ and $v_e$—most apparent in the rapidity of mode exchange through veering or crossing regions—provides a physically grounded indicator that the system is approaching an EP-dominated regime.

In summary, dispersion analysis in non-Hermitian SAFE systems relies critically on distinguishing between analytic structures defined on Riemann surfaces and physically motivated mode tracking. A mode is fundamentally a continuous trajectory on the Riemann surface, and mode tracking aims to reconstruct such continuity from discrete solutions, upon which physical mode labels are assigned.

Accurate dispersion curves therefore require both reliable discrete eigenpair solutions and consistent mode tracking. When these conditions are satisfied, the resulting curves are \emph{spectrally consistent}, meaning they faithfully represent the underlying analytic structure and correctly capture features such as branch topology, modal interactions, and exceptional points. In contrast, inaccuracies in either eigenvalue computation or mode tracking may lead to incorrect branch connectivity, spurious crossings, or misleading curvature, thereby obscuring the true physical spectral structure. Consequently, accurate dispersion curves are essential for obtaining physically meaningful interpretations of wave propagation in non-Hermitian systems.

\section*{CRediT authorship contribution statement}
\textbf{Dong Xiao}: Conceptualization, Methodology, Investigation, Software, Data curation, Formal analysis, Writing - Original draft preparation, Writing - Review Editing, Visualization. \textbf{Zahra Sharif-Khodaei}: Supervision, Writing - Review Editing. \textbf{M. H. Aliabadi}: Supervision, Writing - Review Editing.

\section*{Declaration of competing interest}
The authors declare that they have no known competing financial interests or personal relationships that could have appeared to influence the work reported in this paper.

\section*{Acknowledgements}
The first author acknowledges the financial support from the K. C. Wong Postdoctoral Fellowship, funded by the K. C. Wong Education Foundation.
\section*{Data availability}
The source code and data supporting this study will be made publicly available upon publication at \url{https://github.com/dongxiao96/TopoDisper}. This Python-based tool implements the proposed adaptive homotopy continuation framework, and the repository contains Jupyter notebooks that reproduce all numerical examples presented in this paper.


\setstretch{1.0}
\small{

}

\end{document}